\documentclass[10pt]{article}

\usepackage[T1]{fontenc}
\usepackage{lmodern}
\usepackage{microtype}

\usepackage[margin=1in]{geometry}

\usepackage{amsmath,amsfonts,amssymb,amsthm}

\usepackage{graphicx}
\usepackage{float}
\usepackage{subcaption}
\usepackage{booktabs}
\usepackage[font=small,labelfont=bf]{caption}
\usepackage{diagbox}

\usepackage{xcolor}
\usepackage{verbatim}
\usepackage{url}
\usepackage{soul}
\setulcolor{blue!50!black}

\usepackage[colorlinks=true,
            linkcolor=black,
            citecolor=black,
            urlcolor=black]{hyperref}

\newcommand{\weblink}[2]{\href{#1}{\ul{#2}}}

\theoremstyle{plain}
\newtheorem{theorem}{Theorem}

\theoremstyle{definition}
\newtheorem{definition}{Definition}

\theoremstyle{remark}
\newtheorem{remark}{Remark}

\title{A Fast Nonuniform Solver for the Poisson Equation over a Disk}
\author{
Charlie Pyle\thanks{Department of Mathematics, Texas A\&M University,
College Station, TX, USA.
\href{mailto:charlie.pyle@tamu.edu}{charlie.pyle@tamu.edu}}
\and
Prabir Daripa\thanks{Department of Mathematics, Texas A\&M University,
College Station, TX, USA. \href{mailto:daripa@tamu.edu}{daripa@tamu.edu}}
}
\date{September 2026}

\hypersetup{
  pdftitle={A Fast Nonuniform Solver for the Poisson Equation on a Disk},
  pdfauthor={Charlie Pyle and Prabir Daripa},
  pdfsubject={Numerical analysis; fast Poisson solvers; nonuniform Fourier transforms},
  pdfkeywords={Poisson equation, disk, NUFFT, NUDFT, FFT, Green's function,
               radial recurrence, numerical PDE}
}

\begin{document}
\maketitle

\begin{abstract}
We study fast numerical methods for the Poisson equation on a disk within the FFTRR (Fast Fourier Transform Radial Recurrence) framework, which is built on Green's function representations. Classical FFTRR schemes first apply FFTs in the azimuthal variable and then evaluate mode-by-mode radial recurrences, achieving an overall complexity of \(O(MN\log N)\) on an \(N\times M\) uniform grid, but they require a uniformly spaced azimuthal mesh of $N$ points, limiting their applicability on highly nonuniform sampling patterns. In this work we develop a Nonuniform (NUFFTRR) solver that admits nonuniform grids in both the radial and azimuthal directions while retaining the favorable structure of the original FFTRR formulation. The azimuthal analysis-synthesis step is implemented using either a dense NUDFT least-squares solver or one of two NUFFT-based iterative schemes: a Toeplitz solver using circulant-preconditioned conjugate gradients (PCG), and a preconditioned conjugate gradient for least squares (PCGLS) solver with Pipe–Menon density compensation. These yield azimuthal complexities \(O(N^3 + MN^2)\) for the NUDFT variant and \(O(K_{\mathrm{iter}} MN\log N)\) for the NUFFT-based variants on an \(N\times M\) grid and Krylov iteration count $K_{\mathrm{iter}}$. Numerical experiments on strongly nonuniform meshes demonstrate that the proposed method is robust, spectrally accurate in the azimuthal variable, and competitive in runtime with existing fast Poisson solvers. We make use of vectorization and batched BLAS/GPU-accelerated linear algebra operations to eliminate explicit loops over radii and Fourier modes, while also formulating the azimuthal and radial steps entirely in terms of dense array operations, FFTs, and NUFFTs, allowing for straightforward GPU acceleration. The implementation is released as an open-source Python package at \weblink{https://github.com/CharliePyle4/NUFFTRR_Poisson}{NUFFTRR\_Poisson}, and its methodology can be extended directly to related elliptic problems such as the Helmholtz equation.
\end{abstract}

\section{Introduction}

\paragraph{Background}

The Poisson equation is one of the most widely used elliptic partial
differential equations in scientific computing. It describes the relationship between a scalar field and its underlying source distribution, and it appears in models of gravitation, electrostatics, diffusion, heat conduction, fluid flow, semiconductor devices, and image processing. Its homogeneous counterpart, Laplace's equation, arises in source-free regions and is similarly fundamental to potential theory and steady-state physical models. In physics, the Poisson equation relates gravitational potential to mass density \cite{gravity1Poisson,gravity2Poisson} and governs steady-state heat conduction and diffusion processes \cite{heatTransferPoisson}. It also appears in atmospheric, geophysical, and astrophysical models involving flows on spherical or approximately spherical geometries
\cite{Mei2016Poisson3d,Fornberg1997Poisson3d}. In incompressible fluid
dynamics, a Poisson equation is commonly solved to recover the pressure field from the velocity field during numerical solution of the Navier--Stokes equations \cite{navier1,navier2}. Since this pressure solve is performed at every time step~\cite{Guermondpoisson}, its computational cost can strongly affect the overall cost of a simulation. Poisson equations also occur frequently in computational geometry and computer vision. Poisson surface reconstruction methods recover surfaces from sampled point-cloud or normal-field data \cite{PoissonSurfaceReconstruction}, with applications in three-dimensional scanning, computer-aided design, and image processing. Related Poisson-based formulations are used to extract geometric and shape information from silhouettes and images \cite{silhouettesPoissonModel}. These applications often require repeated solutions of large Poisson problems, motivating the development of accurate and efficient numerical methods. Additional applications arise in electrochemistry, semiconductor modeling, and biological transport. Poisson equations describe electrostatic potentials and steady-state responses to chemical or electrical stimuli in membrane transport problems \cite{electrochemistryPoisson}. They are used in the analysis of MOSFET threshold voltages and in semiconductor device models involving cylindrical geometries, metallic contacts, and dielectric layers
\cite{transistors,SemiconductorNonlinearPoisson}. Fast Poisson solvers have
also been applied to blood-flow simulations in catheterized arteries
\cite{bloodflow}. More generally, radial temperature and density distributions in neutron-absorbing gases can be modeled using Poisson-type equations, with relevance to nuclear-reactor analysis and rocket-engine systems \cite{nuclearandrocketsPoisson}. The breadth of these applications motivates numerical methods that are both computationally efficient and flexible with respect to the underlying grid. In particular, practical measurement and simulation meshes are often not perfectly uniform, creating a need for fast Poisson solvers that can directly
accommodate nonuniform sampling patterns.

\paragraph{Methodology and Related Works}

Fast solvers for the Poisson equation on disks and related circular domains commonly exploit separation of variables in polar coordinates. Fourier transforms are applied in the azimuthal direction, reducing the
two-dimensional problem to a family of one-dimensional radial problems.
When the angular nodes are uniformly spaced, the Fourier analysis and
synthesis steps can be performed efficiently using the FFT, leading to
methods with computational complexity of order $O(MN\log N)$ on an
$N\times M$ polar grid~\cite{FFTRR_Poisson_Disk}. A number of fast disk Poisson solvers combine angular FFTs with finite-difference discretizations in the radial variable. For example, \cite{Lai3orderFDscheme,lyndsey4orderFDscheme} develop third- and fourth-order finite-difference schemes with FFT-based angular transforms. These methods achieve high-order accuracy and favorable computational complexity on uniform polar grids. A principal numerical issue in such discretizations is the coordinate singularity at the origin, which can require shifted radial points or specialized origin treatments. The FFTRR framework provides an alternative disk-based formulation using Green's function representations and radial recurrence relations. Rather than approximating radial differential operators by finite differences, the method computes radial contributions through quadrature and recurrence
relations associated with the Fourier modes. This structure avoids a radial grid shift at the origin and permits the radial integration rule to be replaced by higher-order quadrature rules when additional radial accuracy is required. On a uniform angular grid, the original FFTRR method retains the same FFT-based $O(MN\log N)$ scaling while providing a direct Green's-function treatment of the radial problem~\cite{FFTRR_Poisson_Disk}. Other methods combining Fourier separation with Green's-function-based radial evaluation include the cylindrical-coordinate elliptic solver of Pataki and
Greengard~\cite{PatakiGreengard2011} and the unit-disk Poisson solver developed by Pataki et al.~\cite{PatakiEtAl2013} within a conformal-mapping approach to the Grad--Shafranov equation. The former treats three-dimensional free-space problems and leads to radial modified-Bessel equations, while the latter uses high-order quadrature for the radial Green's-function integrals on the mapped disk. Despite differences in geometry, boundary conditions, and discretization, these methods and the FFTRR formulation share the general principle of combining Fourier analysis with fast evaluation of radial Green's-function integrals. Within the FFTRR framework, related formulations
have also been developed for other elliptic problems and circular geometries \cite{DaripaReviewFFTRR, biharmonic,FFTRR3domains}. In particular, FFTRR-type methods have been applied to both Poisson and Helmholtz equations on interior disks, exterior domains, and annuli \cite{FFTRR3domains}. High-order finite-difference schemes have likewise been extended to disks and spheres with Dirichlet, Neumann, and Robin boundary conditions \cite{Lai2002diskandsphere}. These results indicate that the present nonuniform Fourier extension can potentially be adapted to Helmholtz problems and to other circular or annular domains while retaining the same basic angular-analysis and radial-recurrence structure. The original algorithm and FFTRR formulation for the Poisson problem on the disk was developed in \cite{FFTRR_Poisson_Disk}. We extend this work to nonuniform grids, while also optimizing recurrences and methods. Furthermore, we provide both CPU- and GPU-compatible Python implementations of the code for open-source use through \weblink{https://github.com/CharliePyle4/NUFFTRR_Poisson}{NUFFTRR\_Poisson}.

Although uniform-grid FFT solvers are efficient, they require uniformly spaced angular samples. This assumption can be restrictive when data arise
from irregular measurements, angular encoder errors, adaptive meshes, or intentionally clustered sampling patterns. A common approach is to interpolate nonuniform angular data onto a uniform grid and then apply a conventional FFT solver. However, this preprocessing step introduces an interpolation error that may dominate the total error, particularly when the angular grid is coarse or strongly distorted. The nonuniform discrete Fourier transform (NUDFT) provides a direct alternative by evaluating Fourier expansions at arbitrary nodes. Inverse NUDFT problems can be formulated as dense least-squares systems and solved using dense linear algebra, providing an accurate and straightforward solution. However, their dense transform matrices lead to unfavorable scaling at sufficiently large angular resolutions. The nonuniform fast Fourier transform (NUFFT), first introduced by Dutt and Rokhlin~\cite{DuttRokhlin1993}, used Gaussian-based interpolation and oversampled FFTs to reduce the cost of evaluating nonuniform Fourier transforms. Modern implementations make use of localized spreading, an oversampled FFT, and correction for the interpolation kernel, with libraries such as FINUFFT providing high-accuracy implementations of these transforms while supporting efficient parallel computation \cite{finufft1,finufft2}. Recovering Fourier coefficients from nonuniform samples, however, is generally an inverse problem and commonly requires an iterative method. The present work incorporates these nonuniform Fourier tools into the FFTRR framework. We consider a dense NUDFT least-squares method and two NUFFT-based iterative formulations. Classical circulant-preconditioned Toeplitz least-squares iterations, which solve Toeplitz normal equations using FFT-based preconditioned conjugate gradients (PCG), motivate our Toeplitz-PCG formulation \cite{ChanNagyPlemmons1994}. Our PCGLS variant follows the general conjugate-gradient for least squares (CGLS) and preconditioned conjugate-gradient for least squares (PCGLS) framework for linear least-squares problems and uses Pipe-Menon-style fixed-point updates to construct sampling-density compensation weights for nonuniform Fourier nodes \cite{ArioliGratton2008,PipeMenonDCF}. The choice among these three azimuthal analysis solvers reflects a practical trade-off between grid conditioning, angular resolution, and hardware efficiency. The dense NUDFT achieves low practical overhead on coarse-to-moderate grids via optimized BLAS routines, but scales cubically in $N$. For larger grids, the Toeplitz PCG method exploits the exact Toeplitz structure of the normal operator, together with a circulant preconditioner that is most effective on mildly deformed (e.g., jittered) grids. For severely distorted or clustered geometries where circulant preconditioning breaks down and lumped normal-operator evaluation degrades numerical precision, the PCGLS solver operates directly on paired forward and adjoint NUFFT operators with Pipe-Menon density compensation to ensure numerical stability. A detailed analysis of their conditioning and computational trade-offs is presented in Section~\ref{sec:meshalgs} and Section~\ref{sec:implementation}. All transform and recurrence operations in the present implementation are organized as vectorized array operations, batched dense linear algebra, FFTs, and NUFFTs. This structure avoids explicit loops over Fourier modes and radial nodes where possible, and it is suitable for parallel CPU execution. It also provides GPU acceleration through batched matrix operations and GPU-enabled FFT or NUFFT libraries.

\section{Problem Statement and Method Derivation}

\paragraph{Solution to the Dirichlet Problem}
We first introduce the overarching Dirichlet problem we aim to solve and derive its solution in terms of a particular solution and a homogeneous correction. Extensions to the analogous Neumann problem are discussed in
Section~\ref{sec:neumann}. Below we recall Section 2.1 from Borges \& Daripa~\cite{FFTRR_Poisson_Disk}. Consider the following Poisson equation on a disk with Dirichlet boundary conditions:

\begin{definition}[Poisson Equation with Dirichlet Conditions]
\begin{equation}
\begin{aligned}
\Delta u &= f && \text{in } B, \\
u &= g && \text{on } \partial B.
\end{aligned}
\label{eq:dirichlet-problem}
\end{equation}
where $B = B(0, R) = \{x \in \mathbb{R}^2 : |x| < R\}$.
\end{definition}

We simplify the problem in Eq.~\eqref{eq:dirichlet-problem} by splitting $u$ into a particular solution $v$ and a
homogeneous correction $w$: $v$ absorbs the source term over the unbounded plane, while $w$ supplies the necessary boundary correction. Let $v$ solve
\begin{align*}
\Delta v &= f \quad \text{in } B,
\end{align*}
and let $w$ solve the corresponding homogeneous problem
\begin{align*}
\Delta w &= 0 \quad \text{in } B, \qquad
w = g - v \quad \text{on } \partial B.
\end{align*}
By superposition, the solution of the Dirichlet problem is then
\begin{align*}
u &= v + w.
\end{align*}

The function $v$ can be expressed directly in terms of the source function $f$ using the Green's function for the Laplace operator. The Green's function $G(x,\eta)$ represents the influence at the observation point $x$ caused by the unit source located at $\eta$. In two dimensions, the free-space Green's function for the Laplacian is given by
\begin{align*}
G(x, \eta) &= \frac{1}{2\pi} \log |x - \eta|.
\end{align*}
By superposing the effects from all source points, we obtain
\begin{equation}
\begin{aligned}
v(x) &= \int_{B} f(\eta)\, G(x, \eta)\, d\eta,
\quad x \in B.
\end{aligned}
\label{eq:green-potential}
\end{equation}

A naive evaluation of Eq.~\eqref{eq:green-potential} using quadrature runs into issues due to singularity and would incur $O(N^4)$ complexity on the grid. We, however, make use of Fourier transformations in order to develop a fast and accurate algorithm by separating the numerical scheme into two complementary components. First, we perform an azimuthal Fourier decomposition of the source data \(f\) and boundary data \(g\), mapping samples on each circle \(\rho_\ell\) to Fourier modes \(\widehat{f}_n(\rho_\ell)\) and boundary coefficients \(\widehat{g}_n\) via either FFT or NUFFT-based analysis, and later reconstructing physical values from the solution modes \(\widehat{u}_n(\rho_\ell)\).
Second, for each fixed Fourier mode \(n\) we reduce the original two-dimensional integral for \(v(x)\) to a family of one-dimensional radial integrals, and evaluate these efficiently using closed-form radial recurrences in \(\rho\). The azimuthal Fourier step is described in the next section, followed by the derivation of the radial recurrences and their vectorized implementation.

\section{Azimuthal Component}

In this section, we develop the azimuthal stage of the solver, which decomposes the two-dimensional problem into decoupled radial modes and later reconstructs the physical solution. We work in polar coordinates $(r, \alpha)$ and express the solution and source data in terms of azimuthal Fourier series. Because $u$, $f$, and $g$ are $2\pi$-periodic in $\alpha$, we expand $u$, $f$, and $g$ as
\begin{align*}
u(r,\alpha) &= \sum_{n=-\infty}^{\infty} u_n(r)\,\exp(\mathrm{i} n\alpha), &
u_n(r) &= \frac{1}{2\pi}\int_0^{2\pi} u(r,\alpha)\,\exp(-\mathrm{i} n\alpha)\,d\alpha,\\
f(r,\alpha) &= \sum_{n=-\infty}^{\infty} f_n(r)\,\exp(\mathrm{i} n\alpha), &
f_n(r) &= \frac{1}{2\pi}\int_0^{2\pi} f(r,\alpha)\,\exp(-\mathrm{i} n\alpha)\,d\alpha,\\
g(\alpha) &= \sum_{n=-\infty}^{\infty} g_n\,\exp(\mathrm{i} n\alpha), & 
g_n &= \frac{1}{2\pi}\int_0^{2\pi} g(\alpha)\,\exp(-\mathrm{i} n\alpha)\,d\alpha.
\end{align*}
In the following subsections, we formulate the discrete matrix representation of these transforms and develop fast direct and iterative algorithms to solve the resulting analysis problem on both uniform and irregular angular grids.

\subsection{Matrix Formulation and the Analysis Problem}\label{sec:matrixform}

For each radius $\rho_\ell$ ($\ell = 1, \dots, M$), we collect the azimuthal spatial samples into a vector $f^{(\ell)} \in \mathbb{R}^N$ and the Fourier coefficients into $\widehat{f}^{(\ell)} \in \mathbb{C}^N$. Truncating to the symmetric $N$-mode range $n \in [-N/2, \, N/2-1]$ to match the mode ordering within our implementation, we define the $N\times N$ Fourier matrix
\begin{equation}
\begin{aligned}
A_{k,n} = \exp(\mathrm{i}n\alpha_k), \qquad
k = 0,\dots,N-1,\quad
n = -\frac{N}{2},\dots,\frac{N}{2}-1.
\label{eq:fourier-matrix}
\end{aligned}
\end{equation}
Recovering physical samples from Fourier modes (\emph{synthesis}) and Fourier modes from samples (\emph{analysis}) correspond respectively to
\[
f^{(\ell)} = A\,\widehat{f}^{(\ell)}, \qquad \widehat{f}^{(\ell)} = A^\dagger f^{(\ell)},
\]
where $A^\dagger$ denotes the Moore--Penrose pseudoinverse. Stacking all $M$ radial columns into matrices $F = [f^{(1)}\,\cdots\,f^{(M)}] \in \mathbb{R}^{N\times M}$ and $\widehat{F} = [\widehat{f}^{(1)}\,\cdots\,\widehat{f}^{(M)}] \in \mathbb{C}^{N\times M}$, the collective transformation across all radii becomes
\[
F = A\,\widehat{F}, \qquad \widehat{F} = A^\dagger F.
\]

When the angular nodes are equispaced ($\alpha_k = 2\pi k/N$), the matrix $A$ in Eq.~\eqref{eq:fourier-matrix} is unitary up to scaling \(A^\dagger = A^*/N\). Analysis and synthesis coincide with the classical discrete Fourier transform: applying $A^*/N$ recovers the Fourier coefficients directly. Using standard FFTs along the azimuthal dimension, both analysis and synthesis cost $O(N\log N)$ per radius, giving an overall azimuthal cost of $O(MN\log N)$ across the full disk. However, when the angular nodes $0 \le \alpha_0 < \alpha_1 < \dots < \alpha_{N-1} < 2\pi$ are distinct and shared across all radii, $A$ is a square, nonsingular Vandermonde-like matrix, but is no longer unitary ($A^* \neq A^{-1}$). Synthesis $F = A\widehat{F}$ is evaluated efficiently using a type-2 nonuniform FFT (NUFFT-2). However, applying the adjoint operator $A^*$ via a type-1 NUFFT (NUFFT-1) does \emph{not} recover $\widehat{F}$. Azimuthal analysis thus constitutes a linear inverse problem: for each radial level, we must solve the system $A\,\widehat{f}^{(\ell)} = f^{(\ell)}$. Because the node set $\{\alpha_k\}$ is identical across all radii $\rho_\ell$, the operator $A$ is shared by all $M$ columns of $F$ as well as the boundary vector $g$. For moderate $N$, this inverse problem can be solved directly via dense least squares. For large $N$, we introduce a positive diagonal sampling weight matrix $W = \mathrm{diag}(w_0,\dots,w_{N-1})$ to compensate for angular clustering, define the regularized normal operator
\begin{equation}
\begin{aligned}
T := A^* W A + \lambda I, \qquad (\lambda > 0),
\label{eq:regularized-normal-operator}
\end{aligned}
\end{equation}
and solve the weighted system iteratively using fast transform methods.

\begin{remark}[Nyquist mode storage]
In our implementation, coefficient arrays are stored with length $N+1$ rather than $N$ by explicitly duplicating the Nyquist mode: the entry at $n = -N/2$ is copied to an additional slot at $n = +N/2$, with both endpoint entries halved in amplitude. This symmetric-endpoint convention keeps the stored spectrum symmetric about $n = 0$ for real-valued inputs and simplifies boundary-mode bookkeeping in the radial recurrences of Section~\ref{sec:radialrec}. On an equispaced even-$N$ grid, the endpoint modes alias at the sample nodes ($\exp(-\mathrm{i}N\alpha_k/2) = \exp(\mathrm{i}N\alpha_k/2) = (-1)^k$), making this splitting mathematically exact. On nonuniform nodes, the endpoint exponentials do not generally alias; the length-$(N+1)$ array is therefore an internal algorithmic convention that is exactly inverted prior to synthesis by recombining the halved endpoints ($n = +N/2$ added back into $n = -N/2$) to evaluate the type-2 NUFFT strictly over the original $N$-mode basis $n \in [-N/2, \, N/2-1]$.
\end{remark}

\subsection{Azimuthal Mesh Algorithms}\label{sec:meshalgs}

We now detail the three algorithmic options implemented for the analysis and synthesis stages. These comprise the classical uniform FFT baseline, a direct dense NUDFT solve suited for moderate angular resolutions, and iterative NUFFT formulations designed for large-scale nonuniform meshes. The optimal variant is selected based on the angular resolution, the degree of geometric distortion, and computational considerations.

\paragraph{Uniform Meshes}
On an equispaced grid $\alpha_k = 2\pi k/N$, a single batched forward FFT maps the columns of $F$ and the boundary vector $g$ to their Fourier modes $\widehat{f}_n(\rho_\ell)$ and $\widehat{g}_n$. Following the radial recurrence step, a single batched inverse FFT reconstructs physical solution values $u(\rho_\ell,\alpha_k)$ from $\widehat{u}_n(\rho_\ell)$ across all radii simultaneously, requiring $O(MN\log N)$ total work.

\paragraph{Dense NUDFT on Nonuniform Meshes}
When the angular nodes are irregular but the resolution $N$ is modest, analysis is formulated as the dense least-squares problem
\begin{equation}
\begin{aligned}
\min_{\widehat{f}^{(\ell)}} \bigl\|A \widehat{f}^{(\ell)} - f^{(\ell)}\bigr\|_2,
\qquad \ell = 1,\dots,M.
\label{eq:nudft-least-squares}
\end{aligned}
\end{equation}
Because $A$ is shared across all radii, we stack the data columns and boundary vector into a single joint matrix $[F, \, g] \in \mathbb{R}^{N \times (M+1)}$ and solve the system in a single batched dense linear-algebra call using pivoted QR or SVD. Factoring $A$ costs $O(N^3)$, while applying the factorized inverse to all $M+1$ columns costs $O(MN^2)$, giving an analysis complexity of $O(N^3 + MN^2)$.

After the radial recurrences produce the solution coefficients $\widehat{u}_n(\rho_\ell)$, physical-space values are synthesized on the nonuniform grid via
\[
u(\rho_\ell,\alpha_k) = \sum_{n=-N/2}^{N/2-1} \widehat{u}_n(\rho_\ell)\,\exp(\mathrm{i} n \alpha_k), \qquad k = 0,\dots,N-1,\;\; \ell = 1,\dots,M.
\]
We evaluate this using a single FINUFFT type-2 call with the shared node set $\{\alpha_k\}$ and $M$ simultaneous transforms. For a precision $\varepsilon$, the 1D NUFFT-2 costs $O\bigl(N\log(\varepsilon^{-1}) + N\log N\bigr)$ per transform, yielding an overall synthesis cost of $O\bigl(MN\log(\varepsilon^{-1}) + MN\log N\bigr)$. At moderate resolutions ($N \le 128$), optimized dense BLAS operations run with negligible constant overhead, making the NUDFT variant exceptionally fast in practice despite its asymptotic $N^3$ scaling.

\paragraph{NUFFT+CG on Nonuniform Meshes}
For large angular grids where dense matrix factorization becomes prohibitive, we avoid forming $A$ explicitly and instead recover the Fourier modes iteratively using fast transform operations. We provide two distinct NUFFT-based iterative solvers: a Toeplitz-structured solver using circulant-preconditioned conjugate gradients (PCG) for mildly nonuniform grids, and a preconditioned conjugate gradient for least squares (PCGLS) solver with Pipe-Menon density compensation for more strongly distorted grids.

The normal operator in Eq.~\eqref{eq:regularized-normal-operator} is exactly Toeplitz for arbitrary angular nodes; for mildly perturbed grids, it is particularly well approximated by a circulant preconditioner. Without compensation, local sample clustering over-represents dense angular sectors and severely degrades the conditioning of $A^* A$. To restore approximate discrete orthogonality ($A^* W A \approx I$), we introduce strictly positive density weights $w_j > 0$ that scale inversely with local angular density via an FFT-accelerated periodic kernel density estimate (KDE). The angles $\{\alpha_k\}$ are first binned into a fine uniform grid of $N_{\mathrm{fine}} = \kappa N$ points on $[0,2\pi)$ (oversampling factor $\kappa$), convolved via FFT with a wrapped Gaussian kernel of bandwidth $\sigma=2\pi\beta/N$ (bandwidth parameter $\beta$) to ensure smooth weights, and then interpolated back to the nodes using periodic linear interpolation. The weights are set inversely proportional to this density and normalized so that $\sum_j w_j = 1$.

Because $T$ is Toeplitz in the mode indices with generating vector $t_k = \sum_j w_j \exp(-\mathrm{i} k \alpha_j)$ ($k = -N+1, \dots, N-1$), its action can be evaluated via standard FFT circular convolutions rather than expensive in-loop NUFFTs. We compute this Toeplitz kernel vector $\mathbf{t} = (t_k)$ once up front via a single type-1 NUFFT on the weights and embed it into a $2N \times 2N$ circulant matrix, reducing all subsequent multiplications $T x$ to $2N$-point FFT convolutions plus the regularization $\lambda x$. To accelerate convergence, we construct a positive-definite modification of T.~Chan's optimal circulant preconditioner $\mathcal{P}$ \cite{ChanOptimalCircPrec}. Its first column $c \in \mathbb{C}^N$ averages the diagonals of the unregularized normal operator $A^* W A$:
\[
c_k = \frac{N - k}{N}\, t_k + \frac{k}{N}\, t_{k - N}, \qquad k = 0, \dots, N - 1.
\]
We then form a positive-definite spectral modification by defining
\[
\widetilde{\lambda}_j(\mathcal{P}) = \big| \mathrm{FFT}(c)_j \big| + \mu, \qquad j = 0, \dots, N - 1,
\]
where the shift $\mu > 0$ guarantees strict positive definiteness. The regularization term $\lambda I$ is included in the PCG matrix-vector product but omitted from the circulant kernel.

Because $T$ is identical across all radial rings, all right-hand-side columns $[F, \, g]$ are solved simultaneously in batch. Rather than running independent single-column recurrences, the system is treated as a single stacked block-diagonal operator with the Frobenius inner product $\langle R, Z \rangle_F = \operatorname{Re}\{\operatorname{Tr}(R^* Z)\}$, yielding scalar step sizes $\alpha_k$ and $\beta_k$ shared across the batch. Starting from initial guess $\widehat{F}_0 = 0$, residual $R_0 = B = A^* W [F, \, g]$, preconditioned residual $Z_0 = \mathcal{P}^{-1} R_0$, and search direction $P_0 = Z_0$, each batched PCG iteration executes the following sequence:
\begin{enumerate}
    \item \textbf{Matrix-vector product:} Evaluate $V_k = T P_k$ via $2N$-point FFT circular convolution plus the diagonal shift $\lambda P_k$.
    \item \textbf{Step length:} Compute the step size $\alpha_k = \frac{\langle R_k, Z_k \rangle_F}{\langle P_k, V_k \rangle_F}$.
    \item \textbf{Solution and residual update:} Update the solution $\widehat{F}_{k+1} = \widehat{F}_k + \alpha_k P_k$ and residual $R_{k+1} = R_k - \alpha_k V_k$.
    \item \textbf{Preconditioning:} Apply the circulant preconditioner $Z_{k+1} = \mathcal{P}^{-1} R_{k+1}$ in Fourier space by scaling by $\widetilde{\lambda}_j(\mathcal{P})^{-1}$.
    \item \textbf{Search direction update:} Compute $\beta_k = \frac{\langle R_{k+1}, Z_{k+1} \rangle_F}{\langle R_k, Z_k \rangle_F}$ and update $P_{k+1} = Z_{k+1} + \beta_k P_k$.
\end{enumerate}
Iteration terminates when $\max_\ell (\|R_k^{(\ell)}\|_2 / \|B^{(\ell)}\|_2) \le \tau$ or when $k = K_{\max}$. After $K_{\mathrm{CG}}$ iterations, the total analysis cost is
\[
O\bigl(MN L_{\mathrm{NUFFT}}\bigr) \;+\; O\bigl(K_{\mathrm{CG}}\,MN\log N\bigr).
\]

When the angular mesh exhibits strong clustering, large localized gaps, or systematic geometric distortions, the quality of the circulant preconditioner degrades and PCG convergence slows. Furthermore, evaluating the lumped normal operator $T \approx A^* W A$ can amplify roundoff errors. For these grids, we deploy an unregularized PCGLS solver, which minimizes $\|W^{1/2}(A\widehat{f}^{(\ell)} - f^{(\ell)})\|_2$ directly through alternating forward ($A$, type-2 NUFFT) and adjoint ($A^*$, type-1 NUFFT) transforms. This formulation offers critical numerical stability: the step-size denominator is evaluated directly as the positive weighted norm $\|W^{1/2} q_k\|_2^2$, avoiding loss of positive definiteness or cancellation errors in the lumped operator, and the residual $r_k$ is tracked explicitly in physical space. To accelerate convergence without forming a normal matrix, sampling density compensation is computed up front using the componentwise Pipe-Menon iterative fixed-point updates \cite{PipeMenonDCF}:
\begin{equation}
\begin{aligned}
d_j^{(m)} = \max\Bigl(\operatorname{Re}\bigl\{(AA^* w^{(m)})_j\bigr\}, \, 10^{-12}\Bigr), \qquad
w_j^{(m+1)} = \frac{w_j^{(m)} / d_j^{(m)}}{\sum_{k=0}^{N-1} \bigl(w_k^{(m)} / d_k^{(m)}\bigr)},
\label{eq:pipe-menon-update}
\end{aligned}
\end{equation}
where taking the real part and clamping at $10^{-12}$ guards against floating-point roundoff in $AA^*$, and normalization preserves total unit weight. Dividing by this term dampens oversampled clusters and boosts sparse regions toward discrete orthogonality. Initializing from circular Voronoi trapezoidal intervals provides an accurate initial guess, so that in practice only two fixed-point iterations ($n_{\mathrm{iter}} = 2$) are needed to achieve stable weights. These weights $W = \mathrm{diag}(w_0, \dots, w_{N-1})$ act directly as a spatial preconditioner in the residual adjoint step. All radial columns are solved simultaneously in batch; step sizes $\alpha_k$ are evaluated independently per column, and iteration terminates when the maximum columnwise relative residual satisfies $\max_\ell (\|r_k^{(\ell)}\|_2 / \|r_0^{(\ell)}\|_2) \le \tau$. Starting from $r_0 = f^{(\ell)}$ and $p_0 = s_0 = A^* (W r_0)$, each iteration executes:
\begin{enumerate}
    \item \textbf{Forward transform:} Apply the type-2 NUFFT to obtain $q_k = A p_k$.
    \item \textbf{Step length:} Compute $\alpha_k = \frac{\|s_k\|_2^2}{q_k^* W q_k}$.
    \item \textbf{Solution and residual update:} Update modal coefficients $\widehat{f}_{k+1}^{(\ell)} = \widehat{f}_k^{(\ell)} + \alpha_k p_k$ and spatial residual $r_{k+1} = r_k - \alpha_k q_k$.
    \item \textbf{Adjoint transform:} Apply the type-1 NUFFT to obtain the weighted gradient $s_{k+1} = A^* (W r_{k+1})$.
    \item \textbf{Search direction update:} Compute $\beta_k = \frac{\|s_{k+1}\|_2^2}{\|s_k\|_2^2}$ and update $p_{k+1} = s_{k+1} + \beta_k p_k$.
\end{enumerate}
After $K_{\mathrm{PCGLS}}$ iterations, the total analysis cost is
\[
O\bigl(MN L_{\mathrm{NUFFT}}\bigr) \;+\; O\bigl(K_{\mathrm{PCGLS}}\,MN L_{\mathrm{NUFFT}}\bigr).
\]

For both NUFFT+CG variants, synthesis is identical to the NUDFT case: a single batched type-2 NUFFT reconstructs the solution values $u(\rho_\ell, \alpha_k)$ across all radii in $O\bigl(MN\log(\varepsilon^{-1}) + MN\log N\bigr)$ time.

\section{Radial Recurrences}\label{sec:radialrec}

Having determined the azimuthal Fourier modes $\widehat{f}_n$ and $\widehat{g}_n$, the next stage of the solver computes the radial particular solution and boundary correction for each mode $n$. To achieve optimal complexity, we utilize a fundamental representation theorem from \cite{FFTRR_Poisson_Disk}, which reduces the two-dimensional Green's function integral over the disk into a pair of one-dimensional radial integrals that admit stable closed-form recurrences. We also recall below the steps used there to derive recursive relationships for evaluating one dimensional radial integrals in order to solve our problem efficiently.

\begin{theorem}\label{thm:dirichlet}
The \( n \)th Fourier coefficient $u_n(r)$ of the solution \( u(r, \cdot) \) to the Dirichlet Poisson equation on the disk can be expressed as
\begin{equation}
\begin{aligned}
u_n(r)
&=v_n(r)+\left(\frac{r}{R}\right)^{|n|}
  \bigl(g_n-v_n(R)\bigr),
\qquad 0<r\le R.
\label{eq:dirichlet-modal-solution}
\end{aligned}
\end{equation}
where $g_n$ are the Fourier coefficients of the boundary condition $g$, and the Fourier coefficient of the particular solution $v_n(r)$ is given by
\[
v_n(r) = \int_0^r p_n(r,\rho)\, d\rho + \int_r^R q_n(r,\rho)\, d\rho,
\]
with kernels
\[
p_n(r,\rho) =
\begin{cases}
\rho \log r \, f_0(\rho), & n = 0, \\[6pt]
-\dfrac{\rho}{2|n|} \left(\dfrac{\rho}{r}\right)^{|n|} f_n(\rho), & n \neq 0,
\end{cases}
\qquad
q_n(r,\rho) =
\begin{cases}
\rho \log \rho \, f_0(\rho), & n = 0, \\[6pt]
-\dfrac{\rho}{2|n|} \left(\dfrac{r}{\rho}\right)^{|n|} f_n(\rho), & n \neq 0.
\end{cases}
\]
\end{theorem}
By splitting the radial domain at $\rho = r$, this representation decouples the Green's function into an inner integral over $[0, r]$ and an outer integral over $[r, R]$. In the next subsection, we exploit this separation to derive discrete step-by-step recurrence relations across the radial mesh.

\subsection{Vectorized 1D Radial Recurrences}

We now discretize the disk $B(0, R)$ on an $N \times M$ polar grid with $N$ azimuthal nodes and $M$ distinct (possibly nonuniform) radial coordinates $0 = r_1 < r_2 < \cdots < r_M = R$. Because $v_n(0) = 0$ for every nonzero mode $n$, all nonconstant modal components of the particular solution vanish at the origin; the zeroth mode is treated separately below. For adjacent radial intervals $[r_i, r_j]$ ($r_j > r_i$), we define the local radial increments:
\begin{align*}
C^{i,j}_n &= \int_{r_i}^{r_j} \frac{\rho}{2n} \left( \frac{r_j}{\rho} \right)^n f_n(\rho) \, d\rho, \qquad n < 0, \\
D^{i,j}_n &= - \int_{r_i}^{r_j} \frac{\rho}{2n} \left( \frac{r_i}{\rho} \right)^n f_n(\rho) \, d\rho, \qquad n > 0,
\end{align*}
with the $n = 0$ logarithmic cases defined by
\[
C^{i,j}_0 = \int_{r_i}^{r_j} \rho f_0(\rho)\, d\rho, \qquad
D^{i,j}_0 = \int_{r_i}^{r_j} \rho \log \rho \, f_0(\rho)\, d\rho.
\]

Evaluating the split Green's function integrals of Theorem~\ref{thm:dirichlet} recursively yields an outward sweep for $v_n^-$ and an inward sweep for $v_n^+$: 
\begin{equation}
\begin{aligned}
v_n^{-}(r_1)&=0, &
v_n^{-}(r_j)
&=\left(\frac{r_j}{r_i}\right)^n v_n^{-}(r_i)+C_n^{i,j},
&&n\le 0,\\[4pt]
v_n^{+}(r_M)&=0, &
v_n^{+}(r_i)
&=\left(\frac{r_i}{r_j}\right)^n v_n^{+}(r_j)+D_n^{i,j},
&&n\ge 0.
\label{eq:radial-recurrences}
\end{aligned}
\end{equation}
In the original formulation \cite{FFTRR_Poisson_Disk}, the recurrences in Eq.~\eqref{eq:radial-recurrences} were evaluated sequentially via nested for-loops across radial steps for each Fourier mode. We, however, reformulate the recurrences algebraically to enable vectorized computations using prefix operations. Each sweep is a first-order linear recurrence of the form $y_k = a_k y_{k-1} + C_k$ with $y_0 = 0$, where $a_k = (r_k/r_{k-1})^n$ and $C_k = C_n^{k-1,k}$. Unrolling this recurrence step-by-step:
\begin{align*}
y_1 &= C_1, \\
y_2 &= a_2 C_1 + C_2, \\
y_3 &= a_3 a_2 C_1 + a_3 C_2 + C_3,
\end{align*}
reveals the general closed form:
\[
y_k = \sum_{i=1}^{k} \left(\prod_{j=i+1}^{k} a_j\right) C_i, \qquad \prod_{j=i}^{k} a_j = 1 \text{ for } k < i.
\]
Defining the cumulative product $P_k = \prod_{j=1}^{k} a_j$, the product ratio satisfies $\prod_{j=i+1}^{k} a_j = P_k / P_i$ for $k \ge i$. Factoring $P_k$ outside the summation yields
\begin{equation}
\begin{aligned}
y_k &= P_k\sum_{i=1}^{k}\frac{C_i}{P_i}.
\label{eq:recurrence-prefix-form}
\end{aligned}
\end{equation}
The identity in Eq.~\eqref{eq:recurrence-prefix-form}
converts the sequential radial recurrence into a vectorized
evaluation based on cumulative products and sums (prefix
scans), thereby avoiding explicit loops over the \(M\) radial
nodes. The underlying recurrence \(y_k=a_k y_{k-1}+C_k\) has multipliers \(\lvert a_k\rvert\le 1\), which prevent geometric amplification of previously accumulated perturbations. For the parameter ranges tested, masking prevents the propagation of nonfinite floating-point values. For more extreme modes or more strongly stretched meshes, the original sequential recurrence can be used as a robust fallback. To avoid the coordinate singularity at \(r_1=0\), where \((r_2/r_1)^n\) cannot be evaluated directly, the outward
cumulative product is initialized at the first positive radial
node \(r_2\), with \(P_2^{(n)}\equiv1\):
\[
P_k^{(n)} = \prod_{j=3}^{k} \left(\frac{r_j}{r_{j-1}}\right)^{n}, \qquad k = 2, \dots, M,
\]
which yields the vectorized outward sweep
\begin{equation}
\begin{aligned}
v_n^{-}(r_\ell)
&=P_\ell^{(n)}\sum_{i=2}^{\ell}
  \frac{C_n^{i-1,i}}{P_i^{(n)}},
\qquad n\le 0,\quad \ell=2,\dots,M.
\label{eq:outward-vectorized}
\end{aligned}
\end{equation}
with the boundary condition $v_n^-(r_1) = 0$ enforced analytically. Similarly, for the inward sweep from the outer boundary $r_M = R$, we define the reversed cumulative product
\[
P_k^{(-n)} = \prod_{j=k}^{M-1} \left(\frac{r_j}{r_{j+1}}\right)^{n} = \prod_{j=k}^{M-1} \left(\frac{r_{j+1}}{r_j}\right)^{-n}, \qquad k = 1, \dots, M-1,
\]
where $P_M^{(-n)} \equiv 1$, yielding the vectorized inward sweep
\begin{equation}
\begin{aligned}
v_n^{+}(r_\ell)
&=P_\ell^{(-n)}\sum_{i=\ell}^{M-1}
  \frac{D_n^{i,i+1}}{P_i^{(-n)}},
\qquad n\ge 0,\quad \ell=1,\dots,M-1,
\label{eq:inward-vectorized}
\end{aligned}
\end{equation}
with \(v_n^+(r_M)=0\). Equations~\eqref{eq:outward-vectorized}
and~\eqref{eq:inward-vectorized} replace the two sequential radial sweeps by prefix products and prefix sums. For real-valued source data ($f_{-n} = \overline{f_n}$), the particular solution coefficients $v_n(r_\ell)$ for nonpositive modes are reconstructed for each $\ell = 1, \dots, M$ via:
\begin{equation}
v_n(r_\ell)
=
\begin{cases}
\log r_\ell\,v_0^{-}(r_\ell)+v_0^{+}(r_\ell), & n=0,\\[6pt]
v_n^{-}(r_\ell)+\overline{v_{-n}^{+}(r_\ell)}, & n<0,
\end{cases}
\label{eq:particular-reconstruction}
\end{equation}
where at the origin ($r_1 = 0$) the indeterminate product evaluates to its analytical limit $\lim_{r\to 0^+} \log r \, v^{-}_0(r) = 0$, yielding $v_0(r_1) = v^{+}_0(r_1)$. The positive modes are assigned by algorithmic Hermitian symmetrization, $v_n(r_\ell) = \overline{v_{-n}(r_\ell)}$ ($n > 0$). In our implementation, radial recurrences are evaluated only for nonpositive modes ($n \le 0$), and positive modes are populated via this conjugate relation.

\subsection{Neumann Boundary Conditions}\label{sec:neumann}

The FFTRR framework extends directly to Neumann boundary conditions. The particular solution $v(x)$ and its radial recurrences remain identical; only the harmonic boundary correction $w(x)$ is modified to enforce the prescribed boundary flux. We recall below the theorems and formulas from Borges and Daripa~\cite{FFTRR_Poisson_Disk}.
\begin{definition}[Poisson Equation with Neumann Conditions]
\begin{align*}
\Delta u &= f \quad \text{in } B, \\
\left.\frac{\partial u}{\partial r}\right|_{r=R} &= \psi \quad \text{on } \partial B,
\end{align*}
where $B = B(0, R) = \{x \in \mathbb{R}^2 : |x| < R\}$ and $\psi(\alpha)$ is the prescribed boundary flux.
\end{definition}

Analogous to Theorem~\ref{thm:dirichlet} for Dirichlet conditions, the solution admits the following modal Green's function representation:

\begin{theorem}\label{thm:neumann}
Let $u(r, \alpha)$ solve the Neumann Poisson problem on $B(0, R)$. Then its $n$th azimuthal Fourier mode $u_n(r)$ is given by
\begin{equation}
\begin{aligned}
u_0(r)
&=v_0(r)+\varphi_0,
&&n=0,\\
u_n(r)
&=v_n(r)+\left(\frac{r}{R}\right)^{|n|}
  \left(\frac{R}{|n|}\psi_n+v_n(R)\right),
&&n\ne 0.
\label{eq:neumann-modal-solution}
\end{aligned}
\end{equation}
where $\psi_n$ are the Fourier coefficients of the boundary flux $\psi(\alpha)$, $v_n(r)$ is the particular solution computed via the radial recurrences, and
\[
\varphi_0 = u_0(R) - v_0(R)
\]
is the additive constant determined by a prescribed reference value $u_0(R)$ at the outer boundary.
\end{theorem}

Similar to the radial recurrences, the Neumann boundary update is evaluated without looping over radial nodes or Fourier modes. By casting the wave numbers $|n|$ as a column vector of shape $(N, 1)$ and the radial ratios $(r_\ell/R)$ as a row vector of shape $(1, M)$, the correction term
\[
\left( \frac{r_\ell}{R} \right)^{|n|} \left( \frac{R}{|n|} \psi_n + v_n(R) \right)
\]
is computed across all non-zero modes and all radii simultaneously via two-dimensional array broadcasting, executing in $\mathcal{O}(MN)$ operations on both CPU and GPU.

\begin{remark}[Solvability and Discrete Compatibility]
By the divergence theorem, the Neumann problem admits a solution if and only if the data satisfy the global compatibility condition
\[
\int_B f \, dx = \int_{\partial B} \psi \, ds.
\]
In polar Fourier coordinates, this identity reduces entirely to the $n = 0$ angular mode:
\[
2\pi \int_0^R \rho f_0(\rho)\, d\rho = 2\pi R \psi_0 \quad \Longleftrightarrow \quad C^{1,M}_0 = R\,\psi_0.
\]
Because solutions to the pure Neumann problem are unique only up to an additive constant, fixing the reference value $u_0(R)$ uniquely determines $\varphi_0$ and specifies the solution. In our implementation, rather than rejecting the input or projecting the flux when discrete compatibility is perturbed by quadrature approximation, the user provides the reference value $u_0(R)$; the solver computes $u_0(r) = v_0(r) + (u_0(R) - v_0(R))$, fixing the additive gauge directly.
\end{remark}

\section{Quadrature Schemes}\label{sec:radial_quad}
In this subsection we briefly recall the FFTRR quadrature formulation and discussion within Borges and Daripa~\cite{FFTRR_Poisson_Disk}, and then introduce a vectorized prefix-product implementation and new Simpson-type weights that operate directly on possibly nonuniform radial grids. A straightforward way to approximate the radial integrals is to apply the trapezoidal rule, but this yields only second-order accuracy in the radial step size. One might try to improve the order by inserting auxiliary radii between the original grid points and then using higher-order quadrature formulas to compute $C^{i,i+1}_n$ and $D^{i,i+1}_n$ more accurately. This strategy, however, poses several issues, detailed in \cite{FFTRR_Poisson_Disk}. To avoid these while still improving accuracy, we evaluate the integrals directly on the original mesh using modified two-step recurrences that link consecutive grid points without introducing auxiliary nodes. On a uniform radial grid with spacing $\delta r$ and radii $r_i = (i-1)\,\delta r$, the trapezoidal approximation for $n < 0$ between $r_i$ and $r_{i+1}$ gives
\begin{equation}
\begin{aligned}
C_n^{i,i+1}
&=\frac{(\delta r)^2}{4n}\Bigg[
(i-1)\left(\frac{i-1}{i}\right)^{-n}f_n(r_i)
+i f_n(r_{i+1})\Bigg].
\label{eq:trapezoidal-increment}
\end{aligned}
\end{equation}
with an analogous expression holding for $D^{i,i+1}_n$. In our implementation, the arrays $C_n$ and $D_n$ are computed for all intervals and all modes simultaneously using fully vectorized array operations.

To achieve higher-order accuracy, we evaluate integrals of the form $C^{i-1,i+1}_n$ and $D^{i-1,i+1}_n$ across three consecutive radii $r_{i-1}, r_i, r_{i+1}$, enabling Simpson's rule. For $n < 0$ on a uniform grid, one obtains
\begin{equation}
\begin{aligned}
C_n^{i-1,i+1}
&=\frac{(\delta r)^2}{6n}\Bigg[
(i-2)\left(\frac{i-2}{i}\right)^{-n}f_n(r_{i-1})+4(i-1)\left(\frac{i-1}{i}\right)^{-n}f_n(r_i)
+i f_n(r_{i+1})\Bigg].
\label{eq:simpson-increment}
\end{aligned}
\end{equation}
which yields a higher-order approximation while referencing only existing neighboring circles. Because each increment now spans two radial steps, this leads to modified two-step recurrences. For $n \le 0$ (negative and central modes), we set
\[
v_n^{-}(r_1) = 0,\qquad v_n^{-}(r_2) = C^{1,2}_n,
\]
and for $\ell = 3,\dots,M$,
\[
v_n^{-}(r_\ell) = \left(\frac{r_\ell}{r_{\ell-2}}\right)^{n} v_n^{-}(r_{\ell-2}) + C^{\,\ell-2,\ell}_n.
\]
Similarly, for $n \ge 0$ (positive and central modes), we set
\[
v_n^{+}(r_M) = 0,\qquad v_n^{+}(r_{M-1}) = D^{M-1,M}_n,
\]
and for $\ell = M-2,\dots,1$,
\[
v_n^{+}(r_\ell) = \left(\frac{r_\ell}{r_{\ell+2}}\right)^{n} v_n^{+}(r_{\ell+2}) + D^{\,\ell,\ell+2}_n.
\]

\paragraph{Vectorizing the Simpson Recurrence}
The modified recurrences for Simpson's rule generate a two-step dependence (where $v_n^{-}(r_\ell)$ depends on $v_n^{-}(r_{\ell-2})$), which is not directly compatible with the standard one-step recurrence $y_k = a_k y_{k-1} + C_k$. However, because the recurrence only links indices of the same parity, the radial grid can be algebraically decoupled into independent even ($\ell = 2k$) and odd ($\ell = 2k-1$) sub-sequences:
\begin{align*}
v_n^-(r_{2k}) &= \left(\frac{r_{2k}}{r_{2k-2}}\right)^{n} v_n^-(r_{2k-2}) + C_n^{2k-2,2k}, \qquad (\text{even sub-grid}), \\[4pt]
v_n^-(r_{2k-1}) &= \left(\frac{r_{2k-1}}{r_{2k-3}}\right)^{n} v_n^-(r_{2k-3}) + C_n^{2k-3,2k-1}, \qquad (\text{odd sub-grid}).
\end{align*}
Both branches are now strictly first-order linear recurrences in $k$ and therefore admit the same closed-form prefix-product and prefix-sum factorization derived in Section~\ref{sec:radialrec}. Defining the cumulative geometric ratios over the even and odd sub-grids as
\[
P_k^{(E)} = \prod_{j=2}^k \left(\frac{r_{2j}}{r_{2j-2}}\right)^n, \qquad P_k^{(O)} = \prod_{j=3}^k \left(\frac{r_{2j-1}}{r_{2j-3}}\right)^n \quad \bigl(P_1^{(E)} \equiv 1, \; P_2^{(O)} \equiv 1\bigr),
\]
the explicitly vectorized solutions evaluate in parallel to
\begin{equation}
\begin{aligned}
v_n^-(r_{2k})
&=P_k^{(E)}\left(v_n^-(r_2)
  +\sum_{i=2}^k\frac{C_n^{2i-2,2i}}{P_i^{(E)}}\right),\\[4pt]
v_n^-(r_{2k-1})
&=P_k^{(O)}\sum_{i=2}^k
  \frac{C_n^{2i-3,2i-1}}{P_i^{(O)}}.
\label{eq:simpson-vectorized-outward}
\end{aligned}
\end{equation}
In Eq.~\eqref{eq:simpson-vectorized-outward}, the condition \(v_n^-(r_1)=0\) eliminates the leading constant in the odd subsequence, so
that the first contribution is \(C_n^{1,3}\) at \(r_3\). An identical parity decoupling applies symmetrically to the inward sweep $v_n^+$ ($n \ge 0$), stepping inward from the outer boundary nodes $r_{M-1}$ and $r_M$. In terms of the inward step index $k$, the two branches satisfy
\begin{align*}
v_n^+(r_{M-2k+1}) &= \left(\frac{r_{M-2k+1}}{r_{M-2k+3}}\right)^n v_n^+(r_{M-2k+3}) + D_n^{M-2k+1, M-2k+3}, \\[4pt]
v_n^+(r_{M-2k+2}) &= \left(\frac{r_{M-2k+2}}{r_{M-2k+4}}\right)^n v_n^+(r_{M-2k+4}) + D_n^{M-2k+2, M-2k+4}.
\end{align*}
Defining the inward cumulative geometric ratios in direct analogy with the outward sweep,
\[
P_k^{(-E)} = \prod_{j=2}^k \left(\frac{r_{M-2j+1}}{r_{M-2j+3}}\right)^n, \qquad P_k^{(-O)} = \prod_{j=3}^k \left(\frac{r_{M-2j+2}}{r_{M-2j+4}}\right)^n \quad \bigl(P_1^{(-E)} \equiv 1, \; P_2^{(-O)} \equiv 1\bigr),
\]
the vectorized inward solutions evaluate in parallel to
\begin{equation}
\begin{aligned}
v_n^+(r_{M-2k+1}) &= P_k^{(-E)} \left( v_n^+(r_{M-1}) + \sum_{i=2}^k \frac{D_n^{M-2i+1, M-2i+3}}{P_i^{(-E)}} \right), \\[4pt]
v_n^+(r_{M-2k+2}) &= P_k^{(-O)} \sum_{i=2}^k \frac{D_n^{M-2i+2, M-2i+4}}{P_i^{(-O)}}.
\label{eq:simpson-vectorized-inward}
\end{aligned}
\end{equation}
Here, the outer boundary condition $v_n^+(r_M) = 0$ eliminates the leading constant in the second branch (initializing the sum at $r_{M-2}$ via $D_n^{M-2,M}$), while the first branch is initialized at $r_{M-1}$ from the startup trapezoidal boundary increment $v_n^+(r_{M-1}) = D_n^{M-1,M}$.

\paragraph{Extension to Nonuniform Radial Meshes.}
Because the radial recurrences depend solely on the interval integrals $C_n^{i,j}$ and $D_n^{i,j}$, arbitrary radial spacing $0 = r_1 < r_2 < \dots < r_M = R$ is accommodated directly without altering the recurrence framework. For non-uniform intervals $h_0 = r_i - r_{i-1}$ and $h_1 = r_{i+1} - r_i$ with $H = h_0 + h_1$, we compute the local Simpson quadrature weights analytically:
\begin{equation}\label{eq:nonuniform-three-point-weights}
w_0 = \frac{H(2h_0 - h_1)}{6h_0}, \qquad w_1 = \frac{H^3}{6h_0 h_1}, \qquad w_2 = \frac{H(2h_1 - h_0)}{6h_1}.
\end{equation}
The endpoint weights $w_0$ and $w_2$ are strictly positive provided the local mesh ratio satisfies $1/2 < h_1/h_0 < 2$. For more aggressively graded meshes where this ratio is exceeded, the three-point rule remains algebraically exact for quadratic polynomials, though local weights may become negative. On a uniform radial mesh, this rule reduces to standard Simpson's rule and exhibits fourth-order convergence. On smoothly graded meshes where adjacent step sizes satisfy $h_1 - h_0 = \mathcal{O}(H^2)$, fourth-order global convergence is retained; for arbitrary unequal spacing with $h_1 - h_0 = \mathcal{O}(H)$, the interpolatory quadratic rule is generally third-order globally. Once these localized integrals are formed, the closed-form cumulative-product factorizations apply without modification. This decoupled structure allows radial nodes to be concentrated near boundary layers or origin singularities without introducing a global radial matrix solve or increasing asymptotic complexity.

\paragraph{Higher-Order Schemes and the Coordinate Singularity.}
While we have detailed second- and fourth-order approximations via the trapezoidal and Simpson rules, the flexibility of this integral formulation naturally permits the use of even higher-order quadrature schemes (e.g., Gauss--Legendre panel quadratures or spectral integration). Moreover, finite difference discretizations of the Laplacian in polar coordinates suffer severely from the $1/r$ and $1/r^2$ coordinate singularities at the origin, often requiring grid staggering or pole conditions to maintain stability. In contrast, our integral formulation handles the origin analytically, allowing the boundary conditions at $r_1 = 0$ to be evaluated smoothly and robustly. As a result, higher-order accuracy can be implemented systematically without exacerbating stability issues at the center of the disk.

\section{NUFFTRR Algorithm}

\subsection{Algorithm}\label{sec:algorithm}
We now summarize the full NUFFTRR algorithm for solving the Poisson equation on a disk with arbitrary radial nodes and shared azimuthal sampling (either uniform or nonuniform). The mode- and radius-indexed updates below are stated sequentially for clarity; in the actual implementation, the radial sweeps are evaluated without explicit loops using the cumulative-product prefix-sum factorization of Section~\ref{sec:radialrec} and two-dimensional array broadcasting over $(n,\ell)$. The three azimuthal variants (dense NUDFT, Toeplitz PCG NUFFT, and PCGLS NUFFT) differ only in the analysis inversion used in Step~1.

\textbf{Given:}
\begin{itemize}
  \item Radial nodes \(0 = \rho_1 < \rho_2 < \cdots < \rho_M = R\) (with $\rho_\ell \equiv r_\ell$).
  \item Azimuthal nodes \(\alpha_k \in [0,2\pi)\), \(k = 0,\dots,N-1\), shared across all radii (equispaced \(\alpha_k = 2\pi k/N\) in the uniform case, irregular in the nonuniform case).
  \item Source data values \(f(\rho_\ell,\alpha_k)\) for \(\ell = 1,\dots,M\) and \(k = 0,\dots,N-1\).
  \item Boundary data on \(\partial B\) at \(\rho = R\): Dirichlet boundary values \(g(\alpha_k)\), or Neumann boundary flux \(\psi(\alpha_k)\) along with the reference constant \(u_0(R)\).
\end{itemize}
\noindent
\textbf{Goal:} Compute the physical solution values \(u(\rho_\ell,\alpha_k)\) across the entire disk discretization.

\begin{remark}[Real Data and Algorithmic Hermitian Symmetrization]
For continuous real-valued functions, Fourier coefficients satisfy $f_{-n} = \overline{f_n}$. On an equispaced even-$N$ grid, the two Nyquist exponentials alias at the sampling nodes. On arbitrary nonuniform nodes they do not alias, and coefficients obtained in the asymmetric basis $n = -N/2, \dots, N/2-1$ need not be exactly Hermitian, even for real data. The present real-valued implementation applies an algorithmic Hermitian symmetrization during the radial stage: the recurrences are evaluated for nonpositive modes ($n \le 0$), and the positive interior modes are assigned by
\[
u_n(\rho_\ell) = \overline{u_{-n}(\rho_\ell)}, \qquad 1 \le n < N/2.
\]
The duplicated Nyquist entries are an internal storage convention and are converted to the mode convention required by the synthesis routine before evaluation ($n = +N/2$ added back into $n = -N/2$). Because the single recombined Nyquist exponential $\exp(-\mathrm{i} N \alpha_k / 2)$ is not generally real on nonuniform nodes, the synthesized spatial field $u$ retains a tiny imaginary component and is discarded in real-valued post-processing. This symmetrization is an additional numerical approximation on nonuniform even-$N$ grids. Its effect is small when the highest retained modes are sufficiently resolved. General complex-valued data require independent treatment of both positive and negative modes.
\end{remark}

\medskip
\noindent
\textbf{Step 1: Azimuthal Analysis (Forward Transform).}
Map the data matrix \(F\) and boundary vector \(g\) (or flux vector \(\psi\) for Neumann data) to their azimuthal Fourier coefficients \(\widehat{f}_n(\rho_\ell)\) and \(\widehat{g}_n\) (or \(\widehat{\psi}_n\)), using the analysis method of Section~\ref{sec:meshalgs} appropriate for the angular grid:
\begin{itemize}
  \item \textbf{Uniform mesh:} Apply batched 1D FFTs along the azimuthal axis in $O(MN\log N)$ operations.
  \item \textbf{Nonuniform mesh, dense NUDFT:} Form the shared $N\times N$ Fourier matrix in Eq.~\eqref{eq:fourier-matrix} and solve the least-squares problem in Eq.~\eqref{eq:nudft-least-squares} across all $M+1$ columns simultaneously in $O(N^3+MN^2)$ operations.
  \item \textbf{Nonuniform mesh, Toeplitz PCG NUFFT:} Construct KDE density weights $W$, embed the normal operator from Eq.~\eqref{eq:regularized-normal-operator} into a $2N$ circulant convolution, and solve the resulting system via T.~Chan circulant-preconditioned PCG in $O(MN L_{\mathrm{NUFFT}} + K_{\mathrm{CG}} MN\log N)$ operations.
  \item \textbf{Nonuniform mesh, PCGLS NUFFT:} Compute the Pipe--Menon density weights $W$ using Eq.~\eqref{eq:pipe-menon-update} and apply PCGLS directly to $A$ and $A^*$ via paired forward and adjoint NUFFTs in $O\bigl((2K_{\mathrm{PCGLS}} + 1) MN L_{\mathrm{NUFFT}}\bigr)$ operations.
\end{itemize}

\medskip
\noindent
\textbf{Step 2: Radial increments.}
For trapezoidal quadrature, evaluate single-interval increments $C_n^{i,i+1}$ and $D_n^{i,i+1}$ on each interval $[\rho_i, \rho_{i+1}]$($i = 1,\dots,M-1$) (using Eq.~\eqref{eq:trapezoidal-increment} on uniform grids). For the three-point quadratic rule, evaluate two-step increments $C_n^{i-1,i+1}$ and $D_n^{i-1,i+1}$ across $[\rho_{i-1}, \rho_{i+1}]$ using Eq.~\eqref{eq:simpson-increment} on uniform meshes and the weights in Eq.~\eqref{eq:nonuniform-three-point-weights} on nonuniform meshes.

\medskip
\noindent
\textbf{Step 3: Outward radial sweep.}
For \(n\le0\), evaluate the outward contributions using
\eqref{eq:outward-vectorized}. For Simpson quadrature, use the
parity-decoupled form \eqref{eq:simpson-vectorized-outward}.

\medskip
\noindent
\textbf{Step 4: Inward radial sweep.}
For \(n\ge0\), evaluate the inward contributions using
\eqref{eq:inward-vectorized}. For Simpson quadrature, apply the
reverse-indexed parity decomposition of \eqref{eq:simpson-vectorized-inward}.

\medskip
\noindent
\textbf{Step 5: Combination of particular solutions.}
Reconstruct \(v_n(\rho_\ell)\) for \(\ell=1,\ldots,M\) using
Eq.~\eqref{eq:particular-reconstruction}. The two Nyquist endpoints $n = \pm N/2$ are halved to maintain Hermitian symmetry, and positive modes are recovered via $v_n(\rho_\ell) = \overline{v_{-n}(\rho_\ell)}$ for $n \in [1,\,N/2-1]$.

\medskip
\noindent
\textbf{Step 6: Boundary-condition matching.}
For Dirichlet data, apply Eq.~\eqref{eq:dirichlet-modal-solution} to all modes \(n=-N/2,\ldots,N/2-1\). For Neumann data, apply
Eq.~\eqref{eq:neumann-modal-solution}, with the prescribed reference value \(u_0(R)\) fixing the additive constant.

\medskip
\noindent
\textbf{Step 7: Azimuthal Synthesis (Inverse Transform).}
Reconstruct physical solution values \(u(\rho_\ell,\alpha_k)\) from coefficients \(u_n(\rho_\ell)\). On uniform meshes, this is evaluated via a batched inverse FFT in $O(MN\log N)$ operations. On nonuniform meshes, it is evaluated via a single batched type-2 NUFFT across all $M$ radii in $O(MN L_{\mathrm{NUFFT}})$ operations.

\subsection{Computational Complexity}
We now summarize the computational complexity of each stage. Let $N$ denote the angular node and mode count, $M$ the number of radial nodes, and $K_{\mathrm{NUFFT}}$ the Krylov iteration count ($K_{\mathrm{CG}}$ or $K_{\mathrm{PCGLS}}$). For a requested NUFFT tolerance $\varepsilon$, the transform complexity factor is $L_{\mathrm{NUFFT}} = \log N + \log(\varepsilon^{-1})$, which scales as $\mathcal{O}(\log N)$ for fixed precision.
\begin{itemize}
    \item \textbf{Step 1: Azimuthal analysis.}
    Uniform FFT analysis costs $O(MN\log N)$. Dense NUDFT analysis requires $O(N^3)$ to factor the shared matrix and $O(MN^2)$ to solve for all radial columns, totaling $O(N^3 + MN^2)$. The iterative solvers require $O(MN L_{\mathrm{NUFFT}} + K_{\mathrm{CG}} MN \log N)$ for Toeplitz PCG and $O\bigl((2K_{\mathrm{PCGLS}} + 1) MN L_{\mathrm{NUFFT}}\bigr)$ for PCGLS across all columns, plus a one-time precomputation cost of $O(N L_{\mathrm{NUFFT}})$ for density weights and circulant embeddings.

    \item \textbf{Steps 2--6: Radial operations.}
    Evaluating local quadrature increments, sweeping the prefix-product recurrences, combining modes under Hermitian symmetry, and broadcasting boundary corrections require $\mathcal{O}(1)$ operations per mode per radius. Summed across all $M$ radii and $N$ modes, Steps~2--6 require $\mathcal{O}(MN)$ operations,
    which is strictly lower-order than the angular transform stages.
    \item \textbf{Step 7: Azimuthal synthesis.}
    Uniform FFT synthesis costs $O(MN\log N)$, while nonuniform type-2 NUFFT synthesis costs $O(MN L_{\mathrm{NUFFT}})$.
\end{itemize}
Combining all stages yields the total asymptotic complexities listed in Table~\ref{tab:nufftrr_complexity}. The dense NUDFT solver is optimal for moderate resolutions ($N \le 128$) due to minimal BLAS overhead, while the iterative NUFFT solvers scale with near-linear complexity on large grids.
\begin{table}[!htbp]
\centering
\small
\resizebox{\textwidth}{!}{%
\begin{tabular}{lcccc}
\toprule
& Uniform FFT & Dense NUDFT & Toeplitz PCG & Iterative PCGLS \\
\midrule
Total complexity
& $O(MN\log N)$
& $O\!\left(N^3+MN^2+MN L_{\mathrm{NUFFT}}\right)$
& $O\!\left(MN L_{\mathrm{NUFFT}}
  +K_{\mathrm{CG}}MN\log N\right)$
& $O\!\left((2K_{\mathrm{PCGLS}}+1)
  MN L_{\mathrm{NUFFT}}\right)$ \\
\bottomrule
\end{tabular}%
}
\caption{Total computational complexity of the NUFFTRR algorithm
across azimuthal analysis methods. Here,
$L_{\mathrm{NUFFT}}=\log N+\log(\varepsilon^{-1})$, while
$K_{\mathrm{CG}}$ and $K_{\mathrm{PCGLS}}$ denote the iteration
counts of the respective iterative solvers.}
\label{tab:nufftrr_complexity}
\end{table}

\subsection{Accuracy}

\paragraph{Uniform Case}
On a uniform azimuthal mesh, analysis and synthesis reduce to standard FFT pairs. For smooth data, Fourier coefficients exhibit spectral decay, resolving angular variation down to machine precision with modest $N$; however, under-resolving rapid angular oscillations induces aliasing that limits global accuracy regardless of $M$. Independently, radial errors are governed by the quadrature rules of Section~\ref{sec:radial_quad}: the trapezoidal rule yields quadratic convergence $\mathcal{O}(\delta r^2)$, while the three-point Simpson rule achieves higher-order accuracy. Simpson's rule is inherently fourth-order $\mathcal{O}(\delta r^4)$, but the observed global rate for the lowest modes ($n=0, \pm 1$) can drop to cubic $\mathcal{O}(\delta r^3)$ due to the single trapezoidal startup interval $[r_1, r_2]$ at the origin. For higher Fourier modes ($|n| \ge 2$), the integrand vanishes smoothly at $\rho=0$, and the full fourth-order rate is retained. Accuracy thus reflects a dual trade-off: spectral truncation and aliasing in $N$ versus algebraic quadrature order in $M$.

\paragraph{Nonuniform Case (Azimuthal)}
Nonuniform azimuthal sampling introduces two primary error sources. First, irregular nodes disrupt discrete Fourier orthogonality, elevating $\kappa(A)$ and amplifying numerical sensitivity. Solver selection balances per-iteration cost against geometric robustness: Toeplitz PCG is the fastest iterative variant via $2N$-point FFT convolutions, but its efficiency depends on how accurately the Toeplitz normal operator is approximated by the selected circulant preconditioner, which is most effective on mildly deformed meshes. Conversely, dense NUDFT solves the least-squares system via QR/SVD with backward error governed by $\kappa(A)$ rather than $\kappa(A)^2$, offering superior stability at modest resolutions. For large, severely distorted grids, PCGLS provides a robust alternative by decoupling forward and adjoint NUFFTs, evaluating step lengths via strictly positive sums of squares, and preconditioning the sampling density directly via Pipe-Menon weights. 

Second, from discrete sampling and trigonometric interpolation theory \cite{Landau1967}, stable frequency recovery in finite precision requires sufficiently dense coverage across the circle; local gaps exceeding the nominal Nyquist interval elevate $\kappa(A)$ and severely degrade conditioning regardless of total $N$. Consequently, large angular gaps limit stable frequency recovery relative to a well-resolved uniform grid. This also constrains local angular adaptivity, as refining nodes near localized features while coarsening elsewhere creates destabilizing gaps. Nevertheless, when measurements are fixed on irregular geometries, direct NUFFT and NUDFT solvers are essential to avoid the severe error floors incurred by interpolating data onto uniform grids.

\paragraph{Nonuniform Case (Radial)}
Unlike the azimuthal transform, radial nonuniformity can be leveraged for substantial accuracy gains without introducing a global ill-conditioned system. Because radial discretization relies on local one-dimensional quadrature and decoupled modal recurrences rather than global matrix inversions, redistributing radial nodes carries no global radial-system conditioning cost. For solutions with localized steep gradients such as core concentrations or outer boundary layers, uniform radial meshes under-resolve high-curvature regions and cause global Simpson convergence to degrade. Adapting the radial nodes via power-law clustering ($r_m = R \xi_m^2$), boundary sinh stretching, or Chebyshev--Lobatto clustering equidistributes local truncation error, restoring the optimal fourth-order $\mathcal{O}(M^{-4})$ convergence rate and reducing error by up to three orders of magnitude. Because the analytical quadrature weights and recurrence updates retain the same asymptotic $\mathcal{O}(M)$ radial operation count, these accuracy gains incur negligible computational overhead on both CPU and GPU backends.

\subsection{Implementation}\label{sec:implementation}
The NUFFTRR solver is implemented as an open-source Python package featuring unified CPU and GPU execution pathways dispatched through the top-level \texttt{poisson\_solver} routine via a \texttt{use\_gpu} flag. On CPU, azimuthal FFTs are evaluated using the pyFFTW python wrapper for FFTW (with persistent, thread-safe FFTW plans)\cite{Frigo2005TheDA}, nonuniform transforms via FINUFFT \cite{finufft1,finufft2}, and dense solves via SciPy's \texttt{linalg.lstsq}. The GPU backend mirrors this module-for-module using CuPy's \texttt{cupy.fft} and cuFINUFFT \cite{cufinufft}, with host-device memory transfers handled automatically. For the iterative NUFFT solvers, FINUFFT and cuFINUFFT Guru plans as well as all Krylov work vectors are allocated once during setup and reused across iterations, ensuring that iterations execute entirely via in-place FFTs and array operations. CPU execution utilizes multi-core parallelism via a single \texttt{num\_processors} parameter applied across pyFFTW and FINUFFT, while GPU execution parallelizes across CUDA cores. User configuration is controlled through keyword arguments: \texttt{grid\_type} selects the azimuthal solver (uniform FFT, Toeplitz PCG, or PCGLS), \texttt{use\_nudft\_angular} toggles the dense NUDFT solve,
\texttt{quad\_rule} selects between trapezoidal (1) and Simpson (2)
quadratures, and \texttt{rad\_unif} selects uniform or nonuniform radial
spacing. The iterative solvers expose iteration caps
$K_{\max}$ (\texttt{maxiter\_nufft}), convergence tolerance
$\tau$ (\texttt{tol\_nufft}), transform accuracy
$\epsilon$ (\texttt{eps\_finufft}), Tikhonov regularization
$\lambda$ (\texttt{reg\_param}) for Toeplitz PCG, preconditioner shift
$\mu$ (\texttt{precond\_shift}), and KDE oversampling and bandwidth
parameters $\kappa, \beta$ (\texttt{kde\_oversample},
\texttt{kde\_bandwidth}).

\section{Experiments}

All experiments below were run in a Google Colab environment with an Intel(R) Xeon(R) CPU @ 2.20GHz (6 physical / 12 logical Threads) and 53.0\,GB of system RAM, paired with an NVIDIA L4 GPU (Compute Capability 8.9, 22.0\,GB VRAM) running CUDA 12.9. For each configuration, we report the minimum wall-clock time over five runs, which suppresses startup and caching overhead. In a handful of problems, timings include additional evaluation points beyond those used for the reported accuracy figures; accuracy trends and absolute error values are consistent across these extra points with sufficient regularization and tuning parameters, so we omit these accuracy numbers for brevity. Our first two experiments (Uniform Problems 1 and 2) follow the manufactured solutions and convergence studies of Problems~1 and~5 in Borges \& Daripa~\cite{FFTRR_Poisson_Disk}, providing a direct accuracy baseline and timing on uniform polar grids. The subsequent tests introduce new manufactured solutions and sampling
patterns designed to probe nonuniform angular meshes, radial adaptation, and
CPU/GPU performance of the NUFFTRR variants.

\paragraph{Uniform Problem 1}
In our first example, we use a smooth manufactured solution on the unit circle with a uniform grid:
\begin{align*}
u(x,y) &= 3\exp(x+y)(x-x^2)(y-y^2)+5, \\
f(x,y) &= 6\exp(x+y)xy(-3+x+y+xy).
\end{align*}
This test case provides a baseline for comparing the four solver configurations (Uniform/FFT, Uniform/NUDFT, Uniform/NUFFT-Toeplitz, and Uniform/NUFFT-PCGLS) on a uniform grid, allowing us to verify accuracy and attribute runtimes in terms of transform costs rather than irregular sampling impacts. All tables displaying accuracies are consistent across all algorithms, with the only differences appearing in runtime. 

First, we examine how the relative \(L_2\) error behaves as the angular resolution \(N\) and radial resolution \(M\) are varied independently, using the trapezoidal rule for the radial quadrature and Dirichlet boundary data on the unit circle. We fix the NUFFT parameters to transform accuracy
$\epsilon = 10^{-8}$, CG tolerance $\tau = 10^{-8}$, maximum iterations $K_{\max} = 100$, and regularization $\lambda = 10^{-12}$, although these do not greatly impact runtime or accuracy as the grid is uniform. As shown in Table~\ref{tab:uniform_NM}, the relative \(L_2\) error is essentially independent of \(N\) and shows second order convergence when varying \(M\), confirming that on this smooth problem the uniform angular discretization is already sufficiently accurate, and the dominant source of truncation error is the radial quadrature rather than the angular sampling.

\begin{table}[!htbp]
    \centering
    \begin{tabular}{rccccc}
        \toprule
        N & M = 32 & M = 64 & M = 128 & M = 256 & M = 512 \\
        \midrule
        32  & 1.10e-05 & 2.66e-06 & 6.55e-07 & 1.62e-07 & 4.04e-08 \\
        64  & 1.10e-05 & 2.66e-06 & 6.55e-07 & 1.62e-07 & 4.04e-08 \\
        128 & 1.10e-05 & 2.66e-06 & 6.55e-07 & 1.62e-07 & 4.04e-08 \\
        256 & 1.10e-05 & 2.66e-06 & 6.55e-07 & 1.62e-07 & 4.04e-08 \\
        512 & 1.10e-05 & 2.66e-06 & 6.55e-07 & 1.62e-07 & 4.04e-08 \\
        \bottomrule
    \end{tabular}
    \caption{Uniform Problem 1. Relative \(L_2\) error for Uniform / FFT on the uniform grid as a function of angular resolution \(N\) and radial resolution \(M\).}
    \label{tab:uniform_NM}
\end{table}

\begin{figure}[!htbp]
    \centering
    \includegraphics[width=.75\linewidth]{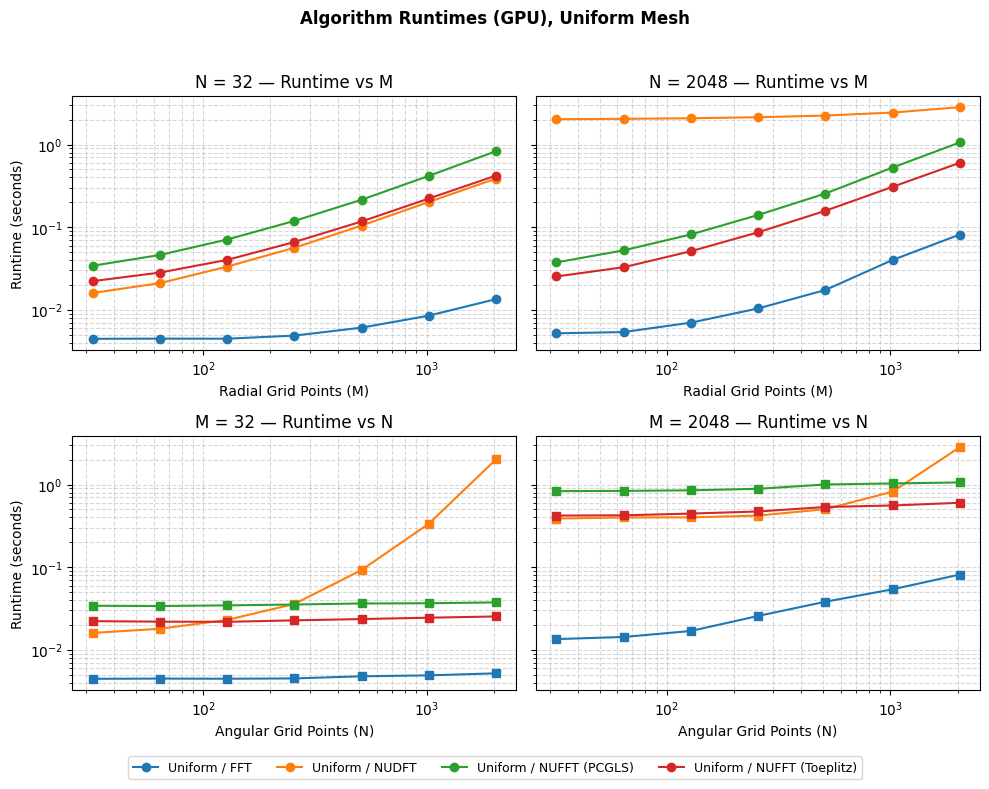}
    \caption{Uniform Problem 1. Runtimes for the four solvers (Uniform/FFT, Uniform/NUDFT, Uniform/NUFFT-Toeplitz, Uniform/NUFFT-PCGLS) on a uniform mesh, shown for the lowest and highest tested values of \(N\) and \(M\).}
    \label{fig:p1runtimes}
\end{figure}

\begin{figure}[!htbp]
    \centering
    \includegraphics[width=1\linewidth]{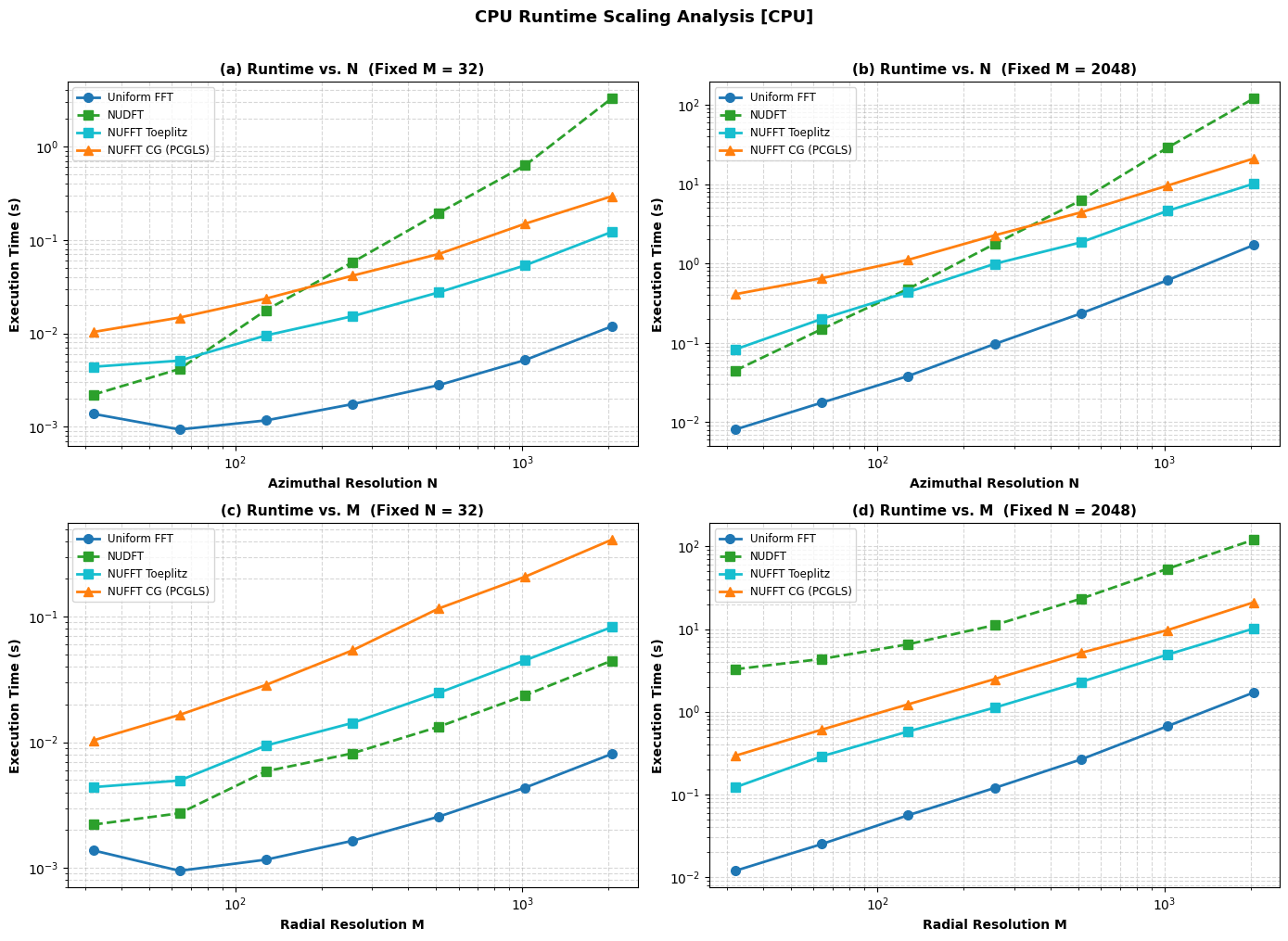}
    \caption{CPU execution time log-log scaling for varying $N$ and $M$.}
    \label{fig:cpu_timing}
\end{figure}

\begin{figure}[!htbp]
    \centering
    \includegraphics[width=1\linewidth]{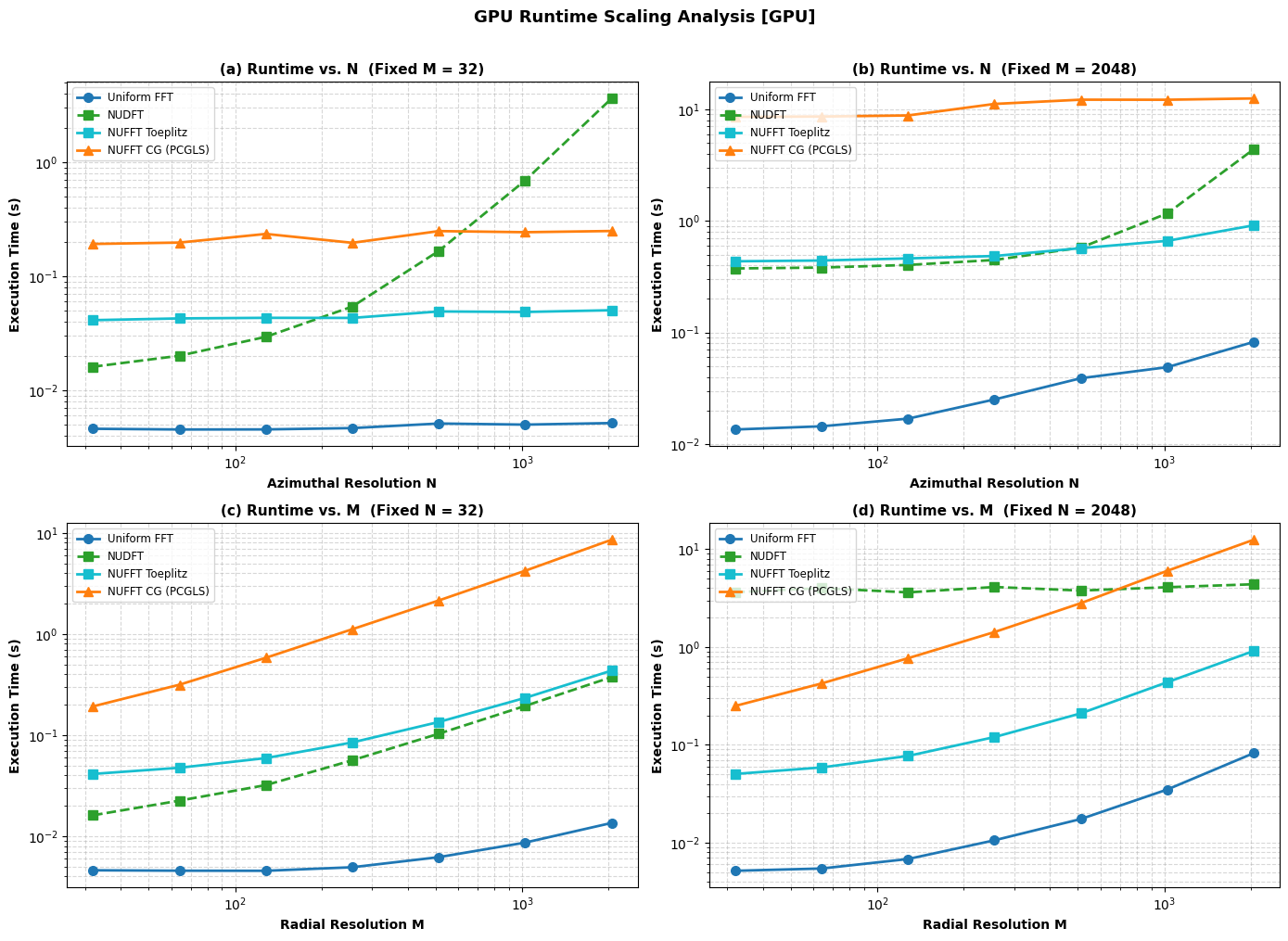}
    \caption{GPU execution time log-log scaling for varying $N$ and $M$.}
    \label{fig:gpu_timing}
\end{figure}

In Figure~\ref{fig:p1runtimes}, we report runtimes at the extreme boundary values of \(N\) and \(M\) across the four algorithms on a GPU backend. The Uniform FFT serves as the optimal baseline and is consistently the fastest solver across all resolutions, while the two NUFFT variants follow identical scaling trends with higher fixed overhead from kernel spreading and iterative Krylov solves. When varying radial resolution \(M\) at a low fixed \(N\), the direct NUDFT outperforms the NUFFT variants across all \(M\). At modest angular sizes, direct matrix-vector evaluation incurs minimal constant overhead and benefits heavily from GPU-accelerated dense matrix operations, whereas the NUFFT methods incur upfront costs-such as kernel spreading, oversampled grid interpolation, and iterative solves-that only amortize when the angular dimension is sufficiently large. At a high fixed \(N\), the NUFFT solvers decisively outperform the direct NUDFT, which rapidly balloons in computational cost. When varying angular resolution \(N\) at a low fixed \(M=32\), runtime curves for the Uniform FFT and NUFFT solvers remain relatively flat at small \(N\) due to fixed GPU kernel launch and dispatch latencies. While the NUDFT is preferable at coarse angular resolutions, increasing \(N\) causes its runtime to quickly surpass the NUFFT variants. A more comprehensive scaling and hardware benchmark across all intermediate grid resolutions is presented later in the timing tests (Figures~\ref{fig:cpu_timing} and~\ref{fig:gpu_timing}).

Next, we fix the angular resolution at \(N=128\) and examine how the choice of quadrature rule and boundary conditions affects the convergence rate as \(M\) increases. As shown in Table~\ref{tab:uniform_quad_bc}, Simpson's rule converges significantly faster than the trapezoidal rule. Although the three-point Simpson rule has a local fourth-order error formula on each subinterval $[r_{i-1}, r_{i+1}]$, the accumulated global error in the radial solver yields cubic convergence ($\mathcal{O}(M^{-3})$), compared to the two-point trapezoidal rule's quadratic convergence. Within the Simpson family, Dirichlet and Neumann boundary conditions produce nearly identical accuracy. For the trapezoidal rule, Neumann errors are approximately an order of magnitude higher because radial quadrature errors in the particular solution $v_n(R)$ enter the Neumann boundary matching formula directly, whereas Dirichlet boundary values pin the solution exactly at $\rho = R$. Regardless, both boundary conditions strictly retain asymptotic second-order convergence.

\begin{table}[!htbp]
    \centering
    \begin{tabular}{rcccc}
        \toprule
        & \multicolumn{2}{c}{Trapezoidal} & \multicolumn{2}{c}{Simpson} \\
        \cmidrule(lr){2-3} \cmidrule(lr){4-5}
        M & Dirichlet & Neumann & Dirichlet & Neumann \\
        \midrule
        32  & 1.10e-05 & 1.88e-04 & 3.14e-06 & 3.18e-06 \\
        64  & 2.66e-06 & 4.55e-05 & 4.02e-07 & 4.05e-07 \\
        128 & 6.55e-07 & 1.12e-05 & 5.08e-08 & 5.10e-08 \\
        256 & 1.62e-07 & 2.78e-06 & 6.39e-09 & 6.40e-09 \\
        512 & 4.04e-08 & 6.91e-07 & 8.00e-10 & 8.01e-10 \\
        \bottomrule
    \end{tabular}
    \caption{Uniform Problem 1. Relative \(L_2\) error for the Uniform/FFT solver as a function of radial resolution \(M\), quadrature rule (trapezoidal vs.\ Simpson), and boundary condition (Dirichlet vs.\ Neumann),  with \(N=128\) fixed.}
    \label{tab:uniform_quad_bc}
\end{table}

\begin{figure}[!htbp]
    \centering
    \includegraphics[width=1\linewidth]{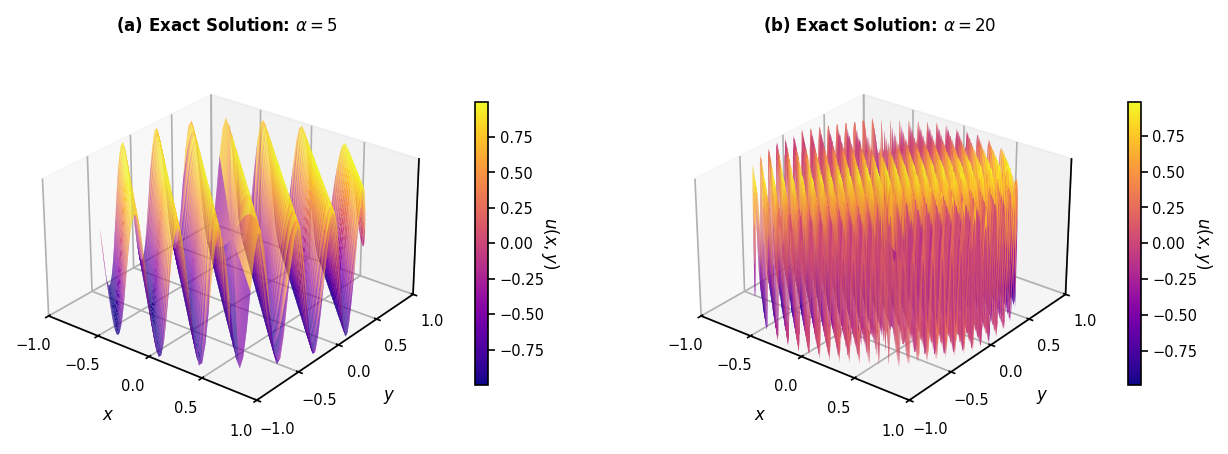}
    \caption{Uniform Problem 2. True solution visualizations for $\alpha=5$ (left) and $\alpha=20$ (right).}
    \label{fig:uniformprob2solution}
\end{figure}

\paragraph{Uniform Problem 2}
We now consider a problem in which the angular resolution plays a more critical role. Consider the following manufactured solution and corresponding source function for the Dirichlet problem:
\begin{align*}
u(x, y) &= \sin(\alpha \pi(x + y)), \\
f(x, y) &= -2 \alpha^2 \pi^2 \sin(\alpha \pi (x+y)).
\end{align*}
Because the functions depend on $\alpha(x+y)$, they oscillate much more rapidly and possess much larger derivatives as $\alpha$ increases. This effect can be seen in Figure~\ref{fig:uniformprob2solution}, where we compare a moderately oscillatory case ($\alpha = 5$) and a highly oscillatory case ($\alpha = 20$).

In Table~\ref{tab:p5_uniform_fft_alpha}, we present the relative \(L_\infty\) errors for the two cases ($\alpha = 5$,  $\alpha=20$), both evaluated on a uniform grid using the trapezoidal rule for the radial integrals.

\begin{table}[!htbp]
    \centering
    \begin{minipage}[t]{0.48\textwidth}
        \centering
        \subcaption{\(\alpha = 5\)}
        \label{tab:p5_uniform_fft_alpha5}
        \resizebox{\linewidth}{!}{%
        \begin{tabular}{rcccccc}
            \toprule
            N & M=64 & M=128 & M=256 & M=512 & M=1024 & M=2048 \\
            \midrule
             64  & 1.4e-02 & 3.4e-03 & 8.4e-04 & 2.1e-04 & 5.2e-05 & 1.4e-05 \\
            128  & 1.4e-02 & 3.4e-03 & 8.4e-04 & 2.1e-04 & 5.2e-05 & 1.3e-05 \\
            256  & 1.4e-02 & 3.4e-03 & 8.4e-04 & 2.1e-04 & 5.2e-05 & 1.3e-05 \\
            512  & 1.4e-02 & 3.4e-03 & 8.4e-04 & 2.1e-04 & 5.2e-05 & 1.3e-05 \\
           1024  & 1.4e-02 & 3.4e-03 & 8.4e-04 & 2.1e-04 & 5.2e-05 & 1.3e-05 \\
           2048  & 1.4e-02 & 3.4e-03 & 8.4e-04 & 2.1e-04 & 5.2e-05 & 1.3e-05 \\
            \bottomrule
        \end{tabular}%
        }
    \end{minipage}%
    \hfill
    \begin{minipage}[t]{0.48\textwidth}
        \centering
        \subcaption{\(\alpha = 20\)}
        \label{tab:p5_uniform_fft_alpha20}
        \resizebox{\linewidth}{!}{%
        \begin{tabular}{rcccccc}
            \toprule
            N & M=64 & M=128 & M=256 & M=512 & M=1024 & M=2048 \\
            \midrule
             64  & 2.6e+01 & 2.5e+01 & 2.5e+01 & 2.5e+01 & 2.5e+01 & 2.5e+01 \\
            128  & 2.2e+00 & 2.1e+00 & 2.1e+00 & 2.1e+00 & 2.1e+00 & 2.1e+00 \\
            256  & 2.8e-01 & 6.5e-02 & 1.6e-02 & 4.0e-03 & 1.0e-03 & 2.5e-04 \\
            512  & 2.8e-01 & 6.5e-02 & 1.6e-02 & 4.0e-03 & 1.0e-03 & 2.5e-04 \\
           1024  & 2.8e-01 & 6.5e-02 & 1.6e-02 & 4.0e-03 & 1.0e-03 & 2.5e-04 \\
           2048  & 2.8e-01 & 6.5e-02 & 1.6e-02 & 4.0e-03 & 1.0e-03 & 2.5e-04 \\
            \bottomrule
        \end{tabular}%
        }
    \end{minipage}
    \caption{Uniform Problem 2. Relative \(L_\infty\) error on the uniform grid for \(\alpha=5\) (left) and \(\alpha=20\) (right).}
    \label{tab:p5_uniform_fft_alpha}
\end{table}
For the smooth case (\(\alpha = 5\)), the error is strictly governed by the radial resolution \(M\). Even at \(N=64\), the Fourier basis fully resolves the azimuthal variations, so increasing \(N\) yields no further accuracy gains. Doubling \(M\) consistently reduces the error by a factor of four, demonstrating the expected second-order quadratic convergence of the radial trapezoidal rule. Conversely, for the highly oscillatory case (\(\alpha = 20\)), coarse angular grids (\(N \le 128\)) under-resolve the rapid angular variations, leading to large errors that cannot be improved by refining \(M\) alone (Table~\ref{tab:p5_uniform_fft_alpha20}). However, once the angular resolution reaches \(N \ge 256\), the Fourier modes fully capture the high-frequency dynamics. At this resolution threshold, angular truncation error ceases to dominate, and clean second-order quadratic convergence in \(M\) resumes across the entire range up to \(M=2048\). This test case highlights that sufficient angular resolution is essential when reconstructing rapidly oscillating functions and validates the stability of the Fourier representation under high-frequency conditions.

\begin{figure}[!htbp]
    \centering
    \includegraphics[width=1\linewidth]{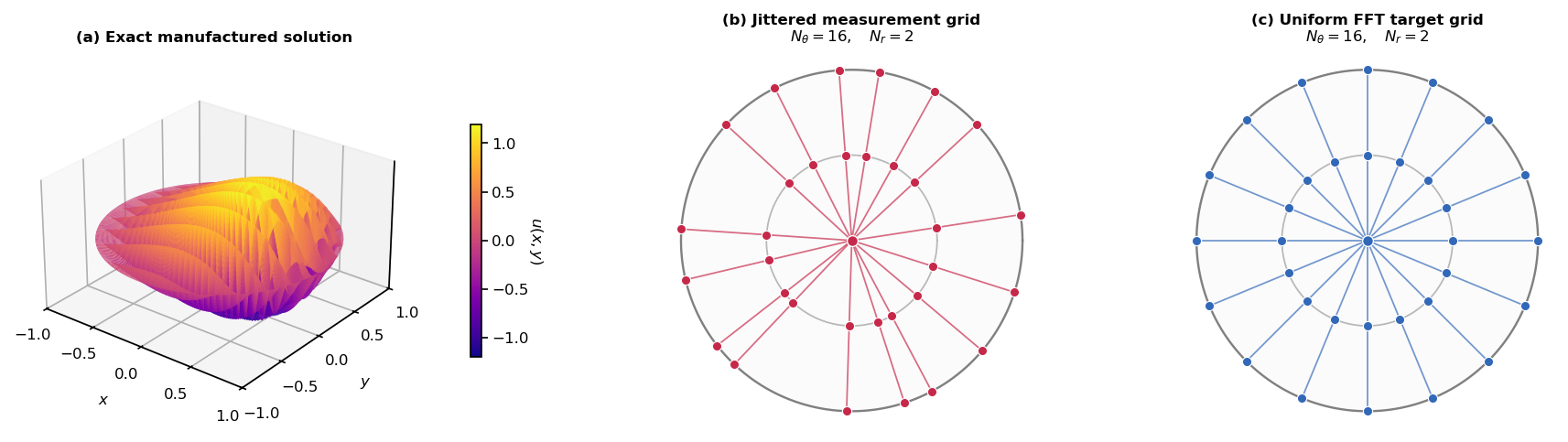}
    \caption{Nonuniform Problem 1. Manufactured solution and angular grids: (a) Exact solution on the unit disk, (b) Jittered polar measurement grid, and (c) Uniform target grid.}
    \label{fig:jittered_solution_and_grids}
\end{figure}

\paragraph{Nonuniform Problem 1 -- Jittered Angular Grid}
To evaluate solver accuracy on an unstructured, quasi-uniform angular grid, we consider the smooth nonseparable manufactured problem on the unit disk ($R=1$) with wave number $k=32$:
\begin{align*}
u(x,y) = (1 - x^2 - y^2) \exp(x) \sin(ky), \qquad u|_{\partial B} = 0.
\end{align*}
The exact forcing $f(x, y) = \Delta u(x, y)$ is obtained analytically as:
\begin{align*}
f(x, y) = \exp(x) \Big[ (1 - k^2)(1 - x^2 - y^2)\sin(ky) - 4(1 + x)\sin(ky) - 4ky \cos(ky) \Big]. 
\end{align*}
Angular nodes are generated by perturbing uniform cell-centered coordinates within each sector:
\begin{align*}
\alpha_j = \frac{2\pi}{N}\left(j + \frac{1}{2}\right) + \delta_j, \qquad |\delta_j| \le 0.45 \left(\frac{2\pi}{N}\right),
\end{align*}
creating an ordered mesh with localized angular compression and dilation. An example grid along with a visualization of the solution is shown in Figure~\ref{fig:jittered_solution_and_grids}.

The NUFFT (Toeplitz) and direct NUDFT solvers operate natively on the jittered angular measurements. For comparison, a Uniform FFT baseline first applies periodic cubic-spline interpolation to resample the jittered data onto an equispaced angular grid before solving. All runs employ radial Simpson quadrature with homogeneous Dirichlet conditions across resolutions $N, M \in \{128, 256, 512, 1024\}$, with NUFFT parameters set to transform accuracy $\epsilon = 10^{-10}$, CG tolerance $\tau = 10^{-10}$,
maximum iterations $K_{\max} = 200$, and regularization $\lambda = 10^{-20}$ (effectively no numerical damping).

\begin{table}[!htbp]
\centering
\small
\begin{tabular}{r cccc cccc}
\toprule
& \multicolumn{4}{c}{NUFFT / NUDFT}
& \multicolumn{4}{c}{Uniform FFT + cubic spline} \\
\cmidrule(lr){2-5} \cmidrule(lr){6-9}
$N \backslash M$
& 128 & 256 & 512 & 1024
& 128 & 256 & 512 & 1024 \\
\midrule
128
& 1.13e-03 & 1.42e-04 & 1.78e-05 & 2.23e-06
& 6.54e-02 & 6.54e-02 & 6.54e-02 & 6.54e-02 \\
256
& 1.13e-03 & 1.42e-04 & 1.78e-05 & 2.22e-06
& 1.83e-03 & 1.51e-03 & 1.52e-03 & 1.52e-03 \\
512
& 1.13e-03 & 1.42e-04 & 1.78e-05 & 2.22e-06
& 1.13e-03 & 1.64e-04 & 8.34e-05 & 8.14e-05 \\
1024
& 1.13e-03 & 1.42e-04 & 1.78e-05 & 2.22e-06
& 1.13e-03 & 1.42e-04 & 1.77e-05 & 2.85e-06 \\
\bottomrule
\end{tabular}
\caption{Nonuniform Problem 1. Relative $L_2$ error on the jittered angular grid across angular ($N$) and radial ($M$) resolutions. The NUFFT and NUDFT methods produce identical errors.}
\label{tab:jittered_grid_l2_error}
\end{table}
As shown in Table~\ref{tab:jittered_grid_l2_error}, NUFFT and NUDFT produce identical accuracy across all resolutions, exhibiting pure radial third-order convergence for $N \ge 128$. In contrast, the Uniform FFT with spline interpolation introduces severe resampling errors that plateau with respect to $M$ on coarse-to-moderate angular grids ($N \le 512$). Only at $N=1024$ does the interpolation gap subside enough to approach the accuracy of the direct nonuniform solvers.

\begin{figure}[!htbp]
    \centering
    \includegraphics[width=1\linewidth]{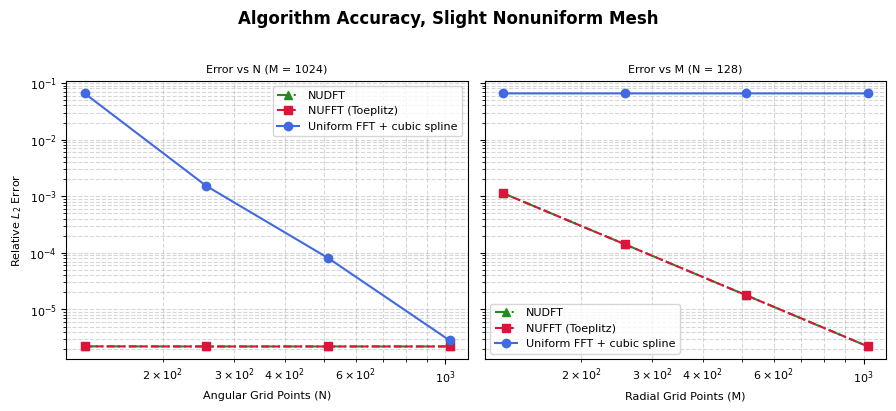}
    \caption{Nonuniform Problem 1. Relative $L_2$ error for the manufactured Poisson problem on a jittered angular grid. Left: error versus the number of angular nodes $N$ at fixed $M=1024$. Right: error versus the number of radial nodes $M$ at fixed $N=128$.}
    \label{fig:jittered_grid_accuracy}
\end{figure}
Figure~\ref{fig:jittered_grid_accuracy} visualizes these decoupled error behaviors. The horizontal plateau in the $N$-refinement plot (left) confirms that the direct solvers are strictly radial-quadrature limited, while the $M$-refinement curves (right) highlight how spline interpolation completely obstructs radial convergence until angular resolution is heavily oversampled.

\begin{figure}[!htbp]
    \centering
    \includegraphics[width=1\linewidth]{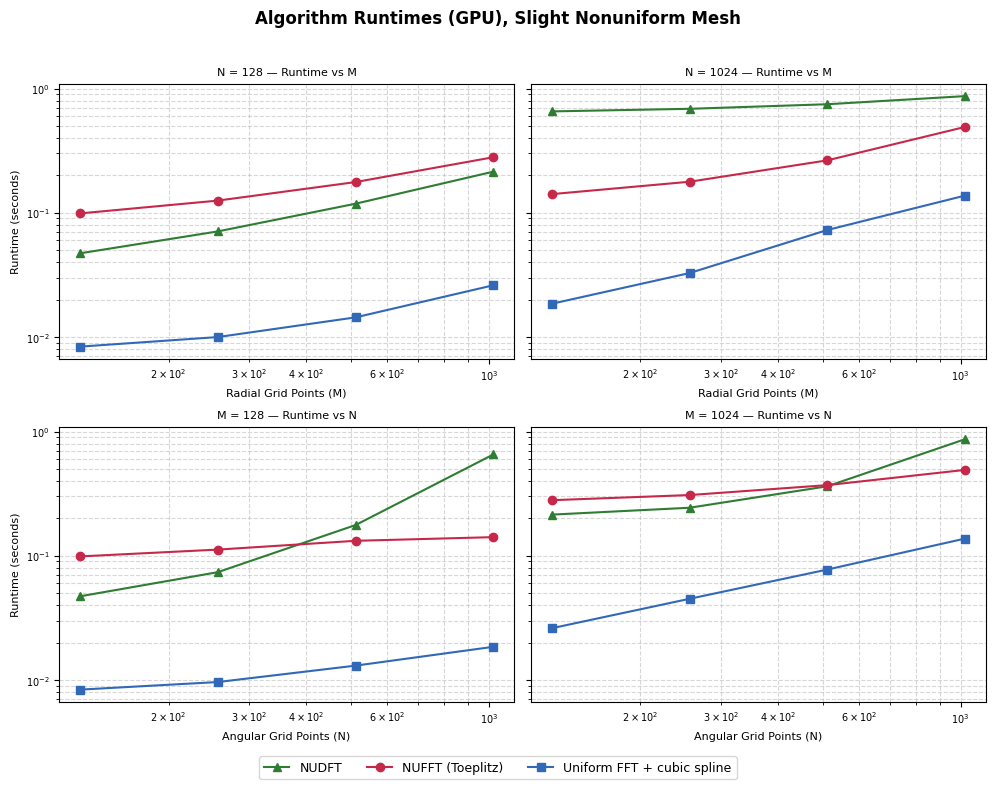}
    \caption{Nonuniform Problem 1. GPU runtime scaling on the jittered angular grid. Top row: runtime versus radial resolution $M$ for fixed $N=128$ and $N=1024$. Bottom row: runtime versus angular resolution $N$ for fixed $M=128$ and $M=1024$.}
    \label{fig:jittered_grid_runtimes}
\end{figure}

Within Figure~\ref{fig:jittered_grid_runtimes} GPU execution times are compared across the three algorithms. The Uniform FFT with periodic cubic-spline interpolation is the fastest overall. For the direct nonuniform solvers, NUDFT is faster than NUFFT Toeplitz at coarse angular resolution ($N=128$) due to minimal setup cost. However, at large angular resolution the NUFFT Toeplitz method displays the best performance, decisively overtaking the quadratic cost of the dense NUDFT ($0.491$\,s vs.~$0.869$\,s at $N=1024, M=1024$). 

Evaluating runtime alongside accuracy reveals a crucial trade-off: at $N=128$ and $M=1024$, both direct NUFFT and NUDFT achieve a relative $L_2$ error of $2.23\times10^{-6}$ in $0.280$\,s and $0.217$\,s, respectively. While the Uniform FFT with cubic spline executes in just $0.022$\,s at this resolution, its accuracy is ruined by interpolation error, yielding a relative error of $6.54\times10^{-2}$. To match the accuracy of the direct nonuniform solvers, the interpolation baseline requires an $8\times$ denser angular measurement budget ($N=1024$), achieving $2.85\times10^{-6}$ in $0.136$\,s. In physical applications where the sensor layout or sampling budget is fixed, increasing angular resolution is often impossible or cost-prohibitive, making direct nonuniform solvers essential for maintaining high accuracy.

\begin{figure}[!htbp]
    \centering
    \includegraphics[width=1\linewidth]{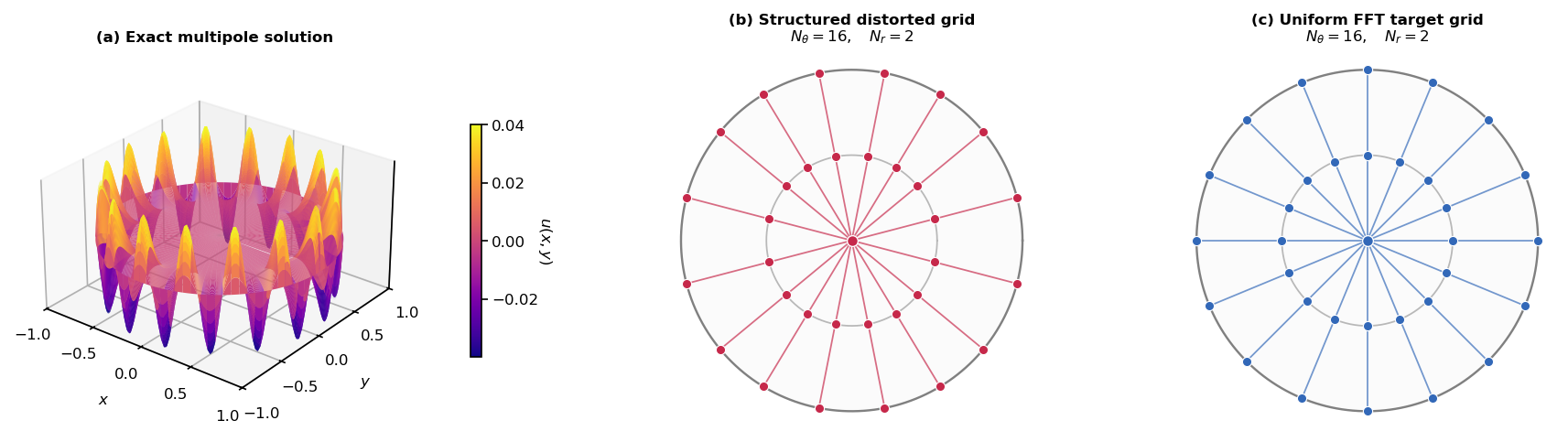}
    \caption{Nonuniform Problem 2. Exact multipole solution and polar angular grids: (a) Exact solution on the unit disk, (b) Structured distorted measurement grid, and (c) Uniform target grid.}
    \label{fig:multipole_solution_and_grids}
\end{figure}

\paragraph{Nonuniform Problem 2 -- Structured Grid}
To evaluate performance on a controlled, deterministic angular distortion, we consider the smooth manufactured single-multipole solution on the unit disk ($R=1$):
\begin{align*}
u(r, \theta) = (1 - r^2)r^m \cos(m(\theta - \theta_0)), \qquad u|_{\partial B} = 0, 
\end{align*}
with exact analytical forcing $f(r, \theta) = \Delta u(r, \theta)$:
\begin{align*}
f(r, \theta) = -4(m + 1) r^m \cos(m(\theta - \theta_0)),
\end{align*}
with mode $m = 16$ and phase $\theta_0 = \pi$. Azimuthal measurement nodes are generated via a structured multipole deformation map:
\begin{align*}
\theta(\xi) = \xi + 0.08 \sin(2\xi) + 0.04 \sin(4\xi), 
\end{align*}
which has a total distortion measure of $0.08(2) + 0.04(4) = 0.32 < 1$. This model represents smooth, non-random angular distortions commonly encountered in physical applications, such as periodic optical lens distortions or angular encoder errors. An example grid along with the solution visualization is shown in Figure~\ref{fig:multipole_solution_and_grids}.

All algorithms are evaluated on the identical distorted measurements across resolutions $N, M \in \{64, 128, 256\}$ using radial Simpson quadrature and Dirichlet boundary conditions. Solver parameters for the iterative NUFFT (PCGLS) are set to transform accuracy $\epsilon = 10^{-10}$, CG tolerance $\tau = 10^{-10}$, and maximum iterations $K_{\max} = 50$. Higher resolutions were omitted because dense, structured multipole clustering introduces ill-conditioning that requires parameter tuning.

\begin{table}[!htbp]
\centering
\begin{tabular}{r ccc ccc}
\toprule
& \multicolumn{3}{c}{NUFFT / NUDFT}
& \multicolumn{3}{c}{Uniform FFT + cubic spline} \\
\cmidrule(lr){2-4} \cmidrule(lr){5-7}
$N \backslash M$
& 64 & 128 & 256
& 64 & 128 & 256 \\
\midrule
64
& 2.77e-05 & 1.75e-06 & 1.09e-07
& 4.19e-02 & 4.19e-02 & 4.19e-02 \\
128
& 2.77e-05 & 1.75e-06 & 1.09e-07
& 1.73e-03 & 1.74e-03 & 1.74e-03 \\
256
& 2.77e-05 & 1.75e-06 & 1.09e-07
& 4.73e-05 & 6.38e-05 & 6.51e-05 \\
\bottomrule
\end{tabular}
\caption{Nonuniform Problem 2. Relative $L_2$ error on the structured distorted grid across angular ($N$) and radial ($M$) resolutions. NUFFT and NUDFT produce identical errors.}
\label{tab:structured_l2error}
\end{table}

As shown in Table~\ref{tab:structured_l2error}, the direct nonuniform methods (NUFFT and NUDFT) yield identical accuracy and substantially outperform the Uniform FFT baseline. For the nonuniform solvers, $N=64$ already fully resolves the azimuthal mode ($m=16$), making the error invariant to further angular refinement and entirely governed by radial Simpson quadrature. In contrast, the Uniform FFT error stalls completely across $M$ at low-to-moderate angular resolutions ($N \le 128$), severely bottlenecked by the cubic-spline interpolation error across distorted spokes. Refining $N$ decreases the interpolation gap, but spline errors continue to limit accuracy even at $N=256$.
To examine practical trade-offs between measurement budget, accuracy, and runtime, Table~\ref{tab:disk_error_comparison} presents execution times and accuracy for a specific problem size on both CPU and GPU backends.

\begin{table}[!htbp]
\centering
\small
\begin{tabular}{lcccccc}
\toprule
Case & $N$ & $M$ & Relative $L_2$ & Relative $L_\infty$ & GPU Runtime (s) & CPU Runtime (s) \\
\midrule
NUFFT (PCGLS) 
    & 64 & 128 & 1.75e-06 & 1.66e-06 & 1.3938 & 0.1076 \\
NUDFT 
    & 64 & 128 & 1.75e-06 & 1.66e-06 & 0.0395 & 0.0125 \\
Uniform FFT + cubic spline 
    & 64 & 128 & 4.19e-02 & 9.37e-02 & 0.0080 & 0.0050 \\
Uniform FFT + cubic spline 
    & 1024 & 128 & 1.58e-06 & 1.62e-06 & 0.0207 & 0.0351 \\
\bottomrule
\end{tabular}
\caption{Nonuniform Problem 2. Relative error and runtime comparison on the structured distorted grid across CPU and GPU backends.}
\label{tab:disk_error_comparison}
\end{table}

At the baseline sampling budget ($N=64$), the Uniform FFT is the fastest method ($0.0050$\,s on CPU, $0.0080$\,s on GPU), but its relative $L_2$ error is more than four orders of magnitude worse ($4.19\times 10^{-2}$) than the direct solvers. To match the accuracy of NUFFT and NUDFT ($1.75\times 10^{-6}$), the Uniform FFT requires a $16\times$ larger angular measurement budget ($N=1024$), taking $0.0351$\,s on CPU and $0.0207$\,s on GPU. In scenarios where measurement density cannot be arbitrarily increased, direct nonuniform methods are required. In this modest resolution regime ($N=64, M=128$), NUDFT is significantly faster than NUFFT (PCGLS) on both CPU ($0.0125$\,s vs.~$0.1076$\,s) and GPU ($0.0395$\,s vs.~$1.3938$\,s). Direct matrix evaluation avoids the iterative Krylov convergence, kernel spreading, and host-device synchronization overhead required by PCGLS. Furthermore, NUDFT executes faster on the CPU than on the GPU ($0.0125$\,s vs.~$0.0395$\,s), as GPU kernel launch latencies limit gains on small grids. Consequently, direct NUDFT is the optimal choice for small-to-moderate distorted grids, while NUFFT's asymptotic advantages become relevant only at much larger scales.

\paragraph{Nonuniform Problem 3 -- Radial Concentration (Origin)}

\begin{figure}[!htbp]
    \centering
    \includegraphics[width=1.0\linewidth]{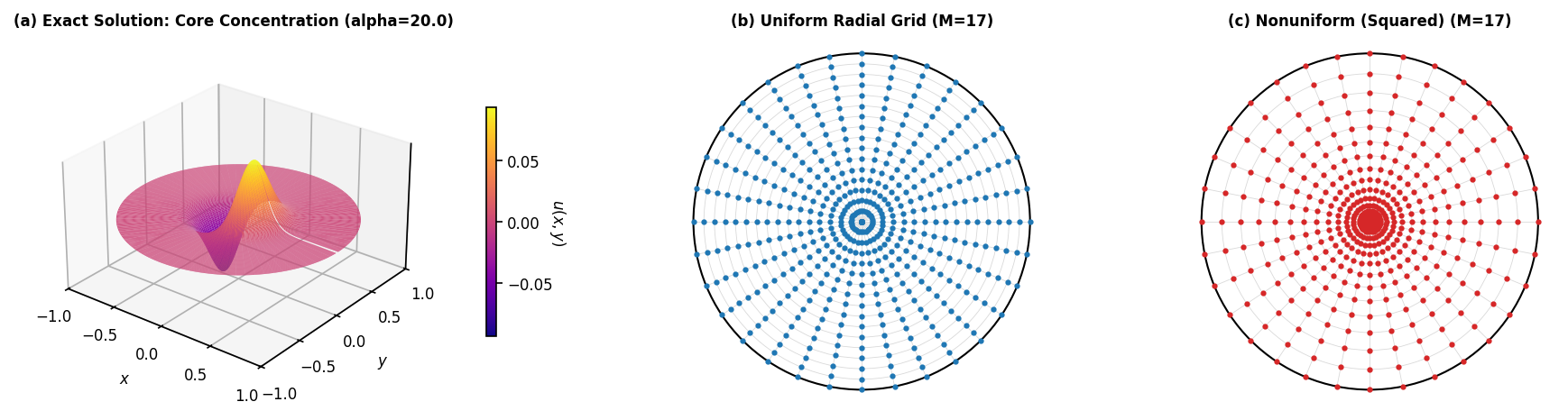}
    \caption{Nonuniform Problem 3. Exact core-concentrated solution and radial grids: (a) 3D surface plot of $u(x,y)$ on the unit disk. (b) Uniform radial grid ($M=17$). (c) Nonuniform (squared) radial grid ($M=17$).}
    \label{fig:p3_solution_and_grids}
\end{figure}

To evaluate the benefits of radial mesh adaptation, we consider a manufactured solution with a localized core peak near the origin:
\[
u(r,\theta) = (R^2-r^2) \exp(-\alpha \frac{r^2}{R^2}) \left(\frac{r}{R}\right)^m \cos\bigl(m(\theta-\theta_0)\bigr),
\]
on the unit disk ($R=1$) with azimuthal mode $m=1$, core concentration parameter $\alpha=20$, and $\theta_0=0$. This solution satisfies homogeneous Dirichlet boundary conditions $u(R,\theta)=0$ and models physical phenomena with intense core excitation, such as focused laser heating or central charge accumulation. The corresponding source term $f(r,\theta) = \Delta u$ is derived analytically:
\[
f(r,\theta) = 4 \exp(-\alpha \frac{r^2}{R^2}) \left(\frac{r}{R}\right)^m \left[ -(m+1)(\alpha+1) + \alpha(\alpha+m+3)\frac{r^2}{R^2} - \alpha^2 \frac{r^4}{R^4} \right] \cos\bigl(m(\theta-\theta_0)\bigr).
\]
Because the solution contains a single azimuthal mode $m=1$, setting $N=32$ resolves the Fourier expansion to machine precision, isolating radial quadrature truncation as the sole error source. While equispaced radial nodes under-resolve steep gradients near the origin, clustering nodes where curvature is highest minimizes local quadrature error. To demonstrate this, we evaluate three radial distributions, visualized alongside the solution in Figure~\ref{fig:p3_solution_and_grids}: an equispaced baseline with uniform nodes $r_j = R \xi_j$; a Chebyshev-Lobatto mesh with dual-boundary clustering $r_j = \frac{R}{2}\bigl(1 - \cos(\pi \xi_j)\bigr)$ near both $r=0$ and $r=R$; and a nonuniform (squared) mesh with power-law core clustering $r_j = R \xi_j^2$ dedicated exclusively to resolving $r \to 0$, where $\xi_j = j/(M-1)$ for $j = 0, \dots, M-1$.

\begin{table}[!htbp]
\centering
\small
\resizebox{\textwidth}{!}{%
\begin{tabular}{r ccc ccc}
\toprule
& \multicolumn{3}{c}{Relative $L_\infty$ Error} & \multicolumn{3}{c}{Relative $L_2$ Error} \\
\cmidrule(lr){2-4} \cmidrule(lr){5-7}
$M$ & Uniform Radial & Chebyshev--Lobatto & Nonuniform (Squared) & Uniform Radial & Chebyshev--Lobatto & Nonuniform (Squared) \\
\midrule
9    & 1.88e-01 & 1.94e-01 & 1.67e-01 & 3.13e-01 & 2.32e-01 & 3.93e-01 \\
17   & 4.88e-02 & 1.55e-02 & 3.65e-03 & 2.42e-02 & 1.78e-02 & 4.47e-03 \\
25   & 1.53e-02 & 2.37e-03 & 6.24e-04 & 5.08e-03 & 2.28e-03 & 7.47e-04 \\
33   & 6.60e-03 & 7.03e-04 & 1.98e-04 & 1.67e-03 & 6.58e-04 & 2.25e-04 \\
49   & 2.00e-03 & 1.29e-04 & 3.79e-05 & 3.48e-04 & 1.23e-04 & 4.31e-05 \\
65   & 8.49e-04 & 4.07e-05 & 1.19e-05 & 1.14e-04 & 3.81e-05 & 1.35e-05 \\
97   & 2.53e-04 & 7.91e-06 & 2.33e-06 & 2.34e-05 & 7.42e-06 & 2.65e-06 \\
129  & 1.07e-04 & 2.49e-06 & 7.37e-07 & 7.63e-06 & 2.34e-06 & 8.35e-07 \\
257  & 1.34e-05 & 1.55e-07 & 4.59e-08 & 5.07e-07 & 1.46e-07 & 5.20e-08 \\
\bottomrule
\end{tabular}%
}
\caption{Nonuniform Problem 3. Relative errors for fixed $N = 32$ with varying radial points $M$, comparing the uniform, Chebyshev--Lobatto, and nonuniform (squared) radial meshes (CPU and GPU produce identical accuracy).}
\label{tab:radial_mesh_comparison}
\end{table}

Table~\ref{tab:radial_mesh_comparison} reports relative $L_\infty$ and $L_2$ errors across radial resolutions $M \in \{9, 17, \dots, 257\}$ at fixed $N=32$. Both nonuniform meshes dramatically outperform the uniform baseline, with the squared mesh achieving the highest precision across all resolutions ($L_\infty \approx 4.59 \times 10^{-8}$ at $M=257$, nearly $300\times$ more accurate than uniform spacing).

\begin{table}[!htbp]
\centering
\small
\resizebox{\textwidth}{!}{%
\begin{tabular}{r ccc ccc}
\toprule
& \multicolumn{3}{c}{CPU Runtime (ms)} & \multicolumn{3}{c}{GPU Runtime (ms)} \\
\cmidrule(lr){2-4} \cmidrule(lr){5-7}
$M$ & Uniform Radial & Chebyshev--Lobatto & Nonuniform (Squared) & Uniform Radial & Chebyshev--Lobatto & Nonuniform (Squared) \\
\midrule
9    & 1.17 & 1.08 & 1.09 & 5.67 & 6.05 & 6.07 \\
17   & 1.18 & 1.14 & 1.13 & 5.84 & 6.15 & 6.18 \\
25   & 1.24 & 1.18 & 1.18 & 5.75 & 6.14 & 6.18 \\
33   & 1.29 & 1.21 & 1.24 & 5.70 & 6.17 & 6.35 \\
49   & 1.36 & 1.29 & 1.30 & 5.69 & 6.21 & 6.17 \\
65   & 1.41 & 1.41 & 1.43 & 5.73 & 6.15 & 6.15 \\
97   & 1.57 & 1.59 & 1.57 & 5.69 & 6.18 & 6.28 \\
129  & 1.75 & 1.70 & 1.73 & 5.71 & 6.18 & 6.31 \\
257  & 2.48 & 2.35 & 2.34 & 6.03 & 6.59 & 6.53 \\
\bottomrule
\end{tabular}%
}
\caption{Nonuniform Problem 3. Solve runtimes (ms) on CPU and GPU for fixed $N = 32$ across radial meshes.}
\label{tab:radial_mesh_runtime}
\end{table}

Table~\ref{tab:radial_mesh_runtime} compares CPU and GPU execution times, where runtimes are virtually identical across all three meshes for each backend. At this small scale CPU execution is faster than the GPU, where kernel launch and stream synchronization latencies dominate. Crucially, nonuniform radial adaptation yields up to three orders-of-magnitude accuracy improvement with zero computational penalty.

\begin{figure}[!htbp]
    \centering
    \includegraphics[width=1\linewidth]{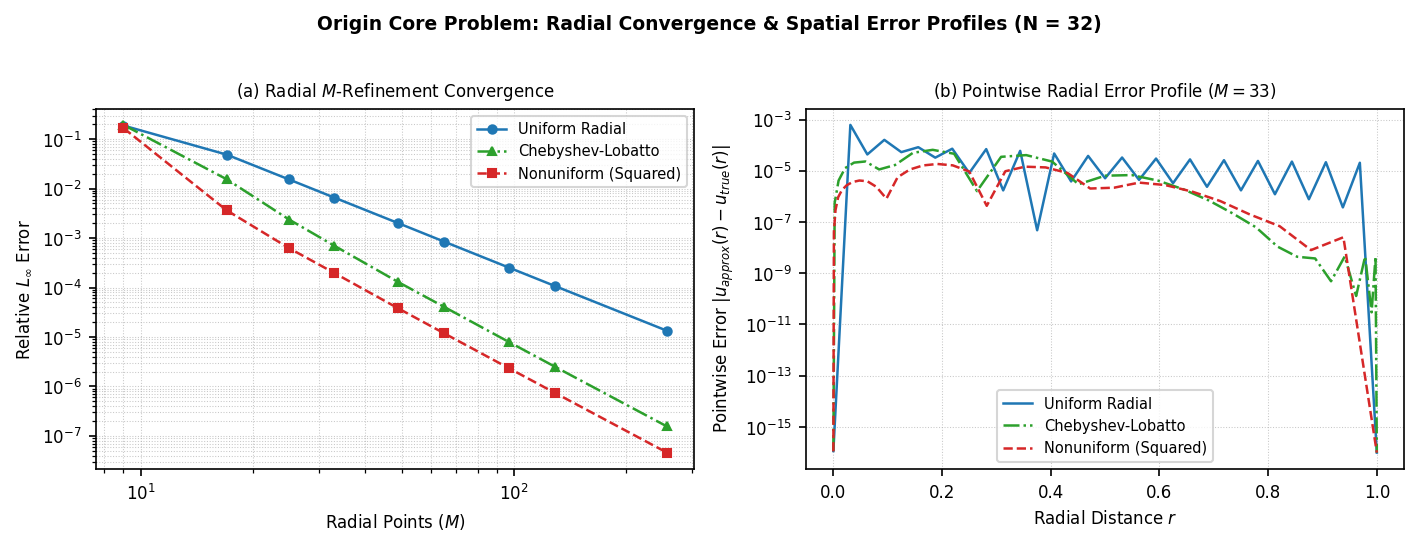}
    \caption{Nonuniform Problem 3. Performance profiles for fixed $N=32$. (a) Log-log $M$-refinement convergence curves. (b) Semilog pointwise radial ray error profiles at $M=33$.}
    \label{fig:p3_convergence}
\end{figure}

Figure~\ref{fig:p3_convergence} illustrates the convergence dynamics and spatial error distributions. The squared and Chebyshev--Lobatto meshes sustain the optimal fourth-order $\mathcal{O}(M^{-4})$ Simpson convergence rate. In contrast, the uniform radial mesh degrades to approximately third-order $\mathcal{O}(M^{-3})$ due to under-resolution of the intense central peak. Pointwise radial error profiles at $M=33$ confirm that core clustering eliminates the severe origin error spike observed on the uniform mesh, establishing a uniformly suppressed error profile across the entire disk radius.

\begin{figure}[!htbp]
    \centering
    \includegraphics[width=1.0\linewidth]{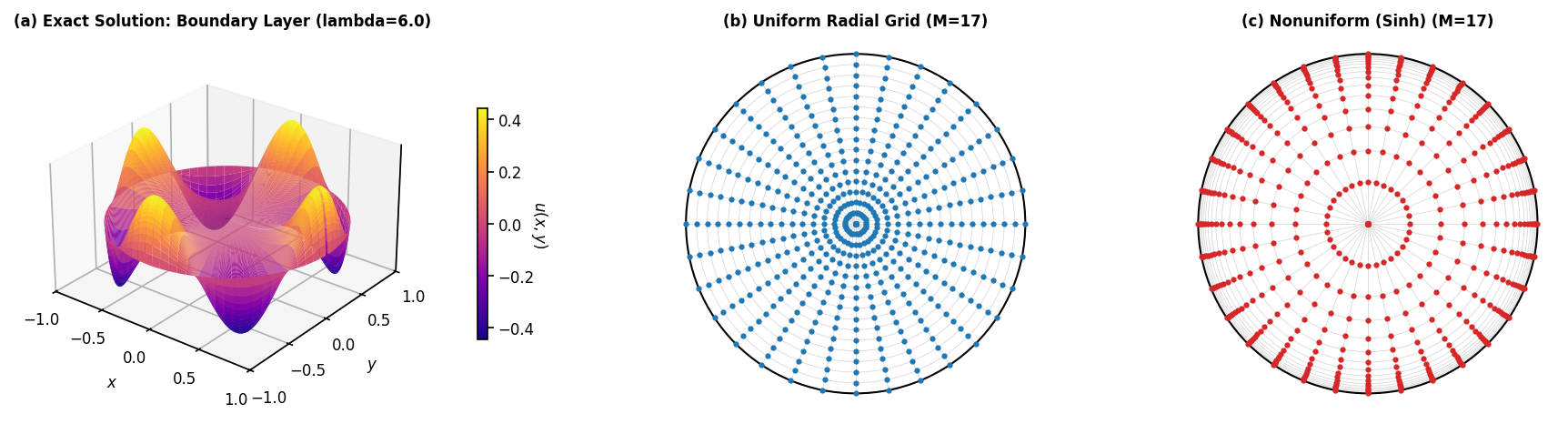}
    \caption{Nonuniform Problem 4. Exact boundary-layer solution and radial grids: (a) 3D surface plot of $u(x,y)$ on the unit disk. (b) Uniform radial grid ($M=17$). (c) Nonuniform (sinh) radial grid ($M=17$).}
    \label{fig:p4_solution_and_grids}
\end{figure}

\paragraph{Nonuniform Problem 4 -- Radial Concentration (Boundary)}
We next examine a problem characterized by a steep outer boundary layer as $r \to R$:
\[
u(r,\theta) = \left(1 - \exp(\beta \left(\frac{r^2}{R^2}-1\right))\right) \left(\frac{r}{R}\right)^m \cos\bigl(m(\theta-\theta_0)\bigr),
\]
on the unit disk ($R=1$) with azimuthal mode $m=4$, boundary layer steepness parameter $\beta=6$, and $\theta_0=0$. The function satisfies homogeneous Dirichlet boundary conditions $u(R,\theta)=0$ and models physical systems with sharp outer gradients, such as thermal boundary layers or skin-effect fields. The corresponding source term $f(r,\theta) = \Delta u$ is derived analytically:
\[
f(r,\theta) = -\frac{4\beta}{R^2} \exp(\beta \left(\frac{r^2}{R^2}-1\right)) \left(\frac{r}{R}\right)^m \left[ m + 1 + \beta \frac{r^2}{R^2} \right] \cos\bigl(m(\theta-\theta_0)\bigr).
\]
As in Problem 3, fixing $N=32$ resolves the $m=4$ Fourier mode to machine precision, isolating truncation errors entirely to radial quadrature. Because equispaced radial nodes severely under-resolve steep boundary gradients near $r \approx R$, clustering nodes where curvature is highest minimizes local quadrature error. To resolve this boundary layer, we evaluate three radial distributions (visualized alongside the solution in Figure~\ref{fig:p4_solution_and_grids}): an equispaced baseline with uniform nodes $r_j = R \xi_j$; a Chebyshev--Lobatto mesh with dual-boundary clustering $r_j = \frac{R}{2}\bigl(1 - \cos(\pi \xi_j)\bigr)$; and a nonuniform (sinh) mesh with hyperbolic sine boundary stretching $r_j = R \bigl( 1 - \frac{\sinh(\gamma (1 - \xi_j))}{\sinh(\gamma)} \bigr)$ with $\gamma = 4.5$ dedicated to resolving $r \to R$, where $\xi_j = j/(M-1)$ for $j = 0, \dots, M-1$.

Table~\ref{tab:radial_mesh_comparison_sinh2} reports relative $L_\infty$ and $L_2$ errors across radial resolutions $M \in \{9, 17, \dots, 257\}$ at fixed $N=32$. Tailoring the stretching parameter $\gamma$ to the boundary layer thickness allows the sinh mesh to achieve the highest accuracy across all resolutions, outperforming the uniform baseline by up to three orders of magnitude ($L_\infty \approx 7.63 \times 10^{-8}$ at $M=257$). The Chebyshev--Lobatto mesh also delivers superior accuracy over the uniform grid by packing nodes at both endpoints.

\begin{table}[!htbp]
\centering
\small
\resizebox{\textwidth}{!}{%
\begin{tabular}{r ccc ccc}
\toprule
& \multicolumn{3}{c}{Relative $L_\infty$ Error} & \multicolumn{3}{c}{Relative $L_2$ Error} \\
\cmidrule(lr){2-4} \cmidrule(lr){5-7}
$M$ & Uniform Radial & Chebyshev--Lobatto & Nonuniform (Sinh) & Uniform Radial & Chebyshev--Lobatto & Nonuniform (Sinh) \\
\midrule
9    & 3.68e-01 & 3.44e-01 & 5.00e-02 & 3.19e-01 & 2.88e-01 & 3.62e-02 \\
17   & 1.30e-01 & 8.94e-03 & 2.51e-03 & 8.59e-02 & 7.86e-03 & 2.08e-03 \\
25   & 5.35e-02 & 1.83e-03 & 4.81e-04 & 3.25e-02 & 1.64e-03 & 3.93e-04 \\
33   & 2.63e-02 & 5.79e-04 & 1.51e-04 & 1.54e-02 & 5.20e-04 & 1.21e-04 \\
49   & 8.94e-03 & 1.14e-04 & 2.96e-05 & 5.08e-03 & 1.02e-04 & 2.28e-05 \\
65   & 4.04e-03 & 3.60e-05 & 9.40e-06 & 2.26e-03 & 3.24e-05 & 6.96e-06 \\
97   & 1.28e-03 & 7.12e-06 & 1.86e-06 & 7.05e-04 & 6.39e-06 & 1.31e-06 \\
129  & 5.59e-04 & 2.25e-06 & 6.02e-07 & 3.05e-04 & 2.02e-06 & 4.08e-07 \\
257  & 7.33e-05 & 1.41e-07 & 7.63e-08 & 3.95e-05 & 1.26e-07 & 3.32e-08 \\
\bottomrule
\end{tabular}%
}
\caption{Nonuniform Problem 4. Relative errors for fixed $N = 32$ with varying radial points $M$, comparing the uniform, Chebyshev--Lobatto, and nonuniform (sinh) radial meshes.}
\label{tab:radial_mesh_comparison_sinh2}
\end{table}

\begin{figure}[!htbp]
    \centering
    \includegraphics[width=1\linewidth]{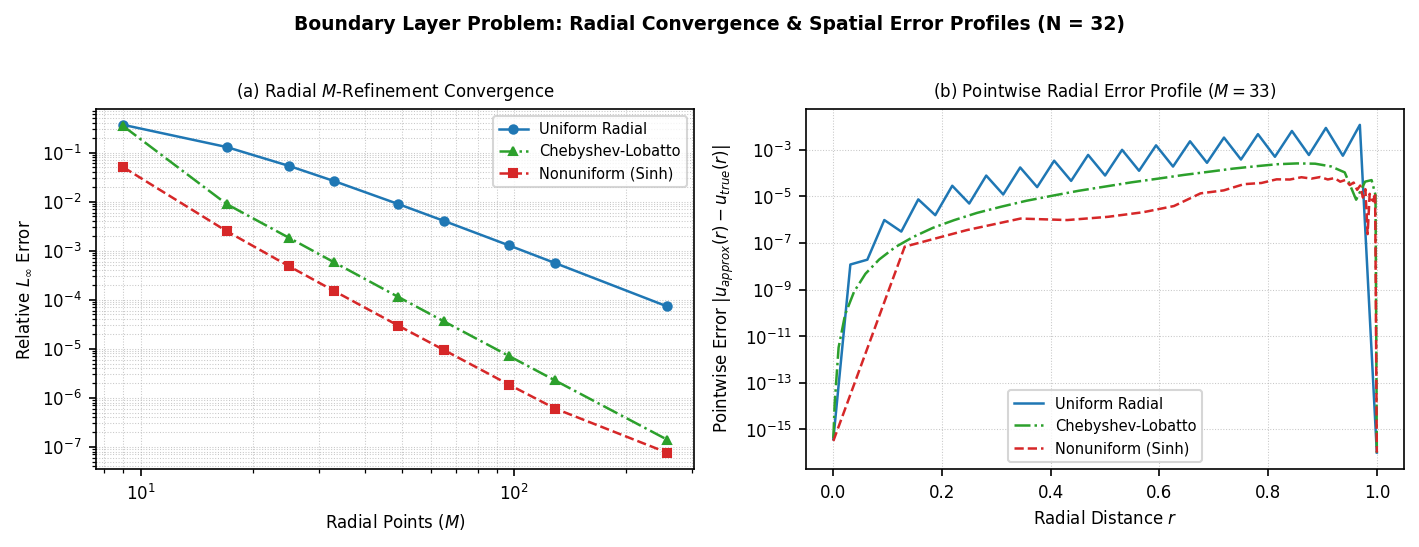}
    \caption{Nonuniform Problem 4. Performance profiles for fixed $N=32$. (a) Log-log $M$-refinement convergence curves. (b) Semilog pointwise radial ray error profiles at $M=33$.}
    \label{fig:p4_convergence}
\end{figure}

Figure~\ref{fig:p4_convergence} illustrates the convergence rates and spatial pointwise error profiles. Both nonuniform meshes sustain optimal fourth-order $\mathcal{O}(M^{-4})$ Simpson convergence slopes. Conversely, the uniform radial mesh degrades to approximately $\mathcal{O}(M^{-2.8})$ convergence because fixed step sizes fail to resolve the steep boundary layer at $r \to 1$. Pointwise radial error profiles at $M=33$ reveal that while all methods enforce Dirichlet conditions at the boundary $r=1$, the uniform mesh develops a severe localized error peak near $r \approx 0.95$. The adapted sinh mesh concentrates nodes near the outer rim, eliminating this boundary error accumulation and maintaining up to two orders of magnitude lower error throughout the radial domain.

\paragraph{CPU and GPU timings}\label{ex:timingtests}

We evaluate runtime scaling across both grid dimensions using the nonseparable manufactured Dirichlet problem of Nonuniform Problem~1 with wave number $k = 2$ on the unit disk $B(0, 1)$. For radial discretization, we employ a uniform radial mesh combined with trapezoidal quadrature. In the angular direction, we evaluate four configurations across an identical resolution grid: the Uniform FFT on an equispaced polar mesh, the direct NUDFT on a jittered nonuniform mesh ($\delta = 0.25$), the fast NUFFT Toeplitz solver on the same jittered mesh, and the iterative NUFFT PCGLS solver on the nonuniform mesh.

For the NUFFT-based solvers, we set the transform accuracy to $\epsilon = 10^{-12}$, Toeplitz Tikhonov regularization to $\lambda = 10^{-10}$, and
KDE preconditioning parameters to oversampling factor $\kappa = 4$, bandwidth factor $\beta = 1.0$, and diagonal shift $\mu = 10^{-3}$. The PCGLS iterations are constrained to a maximum of $K_{\max} = 100$ with stopping tolerance $\tau = 10^{-8}$. The resolution study spans azimuthal dimensions $N \in \{32, 64, 128, 256, 512, 1024, 2048\}$ and radial dimensions $M \in \{32, 64, 128, 256, 512, 1024, 2048\}$. Benchmarking was conducted on both CPU and GPU; CPU computations utilized a single core to obtain clear asymptotic scaling profiles, as multithreading was found to be beneficial only at extreme resolutions while introducing thread synchronization and communication overhead at small-to-moderate scales. It is also important to note these figures represent general behaviors and exact timings will vary under different computer hardware, grid conditioning, or tuning parameters.

In Figure~\ref{fig:cpu_timing}, we show the CPU runtime scaling across both grid dimensions. As expected, the Uniform FFT is the fastest method across all grid sizes. When varying radial resolution $M$, the relative performance depends strongly on the azimuthal resolution $N$. For low fixed $N$, the direct NUDFT consistently outperforms both NUFFT solvers across all $M$. This is because direct matrix multiplication on a small vector avoids the non-negligible setup costs of the NUFFT, such as FINUFFT kernel spreading, oversampled grid interpolation, and iterative Krylov loops. However, for large fixed $N$, the asymptotic scaling of the NUFFT algorithms dominates, making both NUFFT Toeplitz and PCGLS substantially faster than the NUDFT. When varying azimuthal resolution $N$, the NUDFT is faster for small $N$ regardless of the radial scale $M$. However, NUDFT runtimes quickly balloon as $N$ grows, becoming significantly slower than the NUFFT methods. Between the two NUFFT variants, NUFFT Toeplitz is consistently faster than NUFFT PCGLS because circulant embedding solves the angular problem directly via FFTs, bypassing the repeated inner products and iterative evaluations required by PCGLS.

In Figure~\ref{fig:gpu_timing}, we show the corresponding GPU runtime scaling. As on the CPU, the Uniform FFT remains the fastest solver overall, but the behavior of the nonuniform solvers reveals key architectural differences from the CPU. When varying radial resolution $M$, the NUDFT is the optimal choice across all $M$ when $N$ is small, as dense matrix-vector operations are heavily accelerated by GPU dense linear algebra. For high fixed $N$, NUDFT is slower than NUFFT Toeplitz at small $M$; however, as $M$ increases, the NUDFT runtime curve appears remarkably flat. This flatness occurs because batching operations across radial rings $M$ fully saturates GPU memory bandwidth and compute cores without incurring a proportional runtime penalty, whereas NUFFT and PCGLS exhibit steeper growth with $M$. When varying azimuthal resolution $N$, the NUDFT is initially the fastest nonuniform method at small $N$, but its runtime explodes at larger $N$ due to quadratic scaling. Across moderate-to-large $N$, NUFFT Toeplitz is the fastest nonuniform solver on the GPU. Meanwhile, NUFFT PCGLS suffers on the GPU relative to Toeplitz because sequential Krylov iterations require frequent host-device synchronizations and multiple kernel launches per iteration, limiting GPU hardware occupancy.

\paragraph{Other testing}
Beyond the accuracy studies in the present work, the original FFTRR paper \cite{FFTRR_Poisson_Disk} tested five additional manufactured problems, namely Problems 2, 3, 4, 6 and 7 there, that probe complementary aspects of the algorithm and are not repeated here. The same problems were solved using the new  implementation and code proposed in this paper. We obtained the same results, so we briefly summarize the findings. Problem 2 used a solution with a discontinuity in its ``2.5'' derivative and still recovered the expected trapezoidal and Simpson convergence rates in $M$, showing that the radial recurrence is robust to mild smoothness loss in the data as long as the Fourier truncation itself remains adequate. Problems 3 and 4 solved the same class of smooth, non-symmetric solutions on disks of varying radius ($R = 0.5, 1, 2$) at fixed $N$, and found that accuracy improves for smaller $R$ (denser effective point spacing) and degrades for larger $R$ even when $N$ is doubled, highlighting that the domain radius directly rescales the effective radial resolution and must be accounted for when choosing $M$ and $N$ jointly. Problem 6 examined a sharply peaked, near-Gaussian solution and showed that an under-resolved azimuthal grid produces persistent aliasing errors localized near the peak regardless of how much the radial grid $M$ is refined, and that this aliasing only disappears once $N$ is increased enough to resolve the peak's angular content, after which convergence becomes global; this is the same qualitative mechanism seen in our experiment within Uniform Problem 2, but here it is driven by spatial localization rather than a globally rapid azimuthal mode. Problem 7 imposed discontinuous Dirichlet boundary data and showed that the resulting pointwise error stays confined near the discontinuities in $\alpha$ even as $N$ grows, while the error measured away from these points keeps converging, indicating that the method degrades only locally in the presence of non-smooth boundary conditions.

\section{Conclusion}
We have presented NUFFTRR, a fast Green's-function-based solver for the Poisson equation on a disk that accommodates nonuniform discretizations in both the azimuthal and radial coordinates. By combining nonuniform Fourier analysis with closed-form, fully vectorized radial recurrences evaluated via cumulative array operations, the method eliminates explicit mode-by-mode loops and maps seamlessly to batched GPU execution. Because the radial recurrences decouple across Fourier modes and rely strictly on localized one-dimensional quadrature, tailoring radial nodes to resolve concentrated peaks near the origin or steep boundary layers improves accuracy by up to three orders of magnitude with no global radial solve and negligible computational overhead. Across azimuthal discretizations, dense NUDFT provides optimal execution times on coarse-to-moderate grids ($N \le 128$) due to low BLAS overhead, while iterative NUFFT formulations decisively surpass NUDFT at larger resolutions ($N \ge 256$), with Toeplitz PCG providing rapid circular convolution on jittered meshes and PCGLS supporting convergence under more structured clustering. Because this operator-decoupling structure is largely independent of the specific boundary geometry, the NUFFTRR framework extends naturally to annular and exterior domains as well as to related elliptic systems such as the Helmholtz equation.

\bibliographystyle{plain}
\bibliography{references}

\end{document}